
%
%
\def\a{\alpha}

\def\ad{{\rm ad}}
\def\b{\beta}

\def\brel{\buildrel}
\def\bsk{\bigskip}
\def\bu{\bullet}

\def\CC{{\bf C}}
\def\cA{{\cal A}}

\def\cC{{\cal C}}

\def\cE{{\cal E}}
\def\cF{{\cal F}}
\def\cG{{\cal G}}
\def\cH{{\cal H}}

\def\cK{{\cal K}}
\def\cl{\colon}

\def\cM{{\cal M}}

\def\cO{{\cal O}}
\def\cod{{\rm cod}}
\def\cP{{\cal P}}
\def\cQ{{\cal Q}}

\def\cS{{\cal S}}
\def\cU{{\cal U}}

\def\cW{{\cal W}}
\def\cX{{\cal X}}

\def\cZ{{\cal Z}}

\def\D{\Delta}

\def\deg{{\rm deg}}
\def\Def{{\rm Def}}
\def\det{{\rm det}}

\def\e{\epsilon}

\def\Ext{{\rm Ext}}
\def\es{\emptyset}

\def\g{\gamma}
\def\G{\Gamma}

\def\hb{\hbox}
\def\Hom{{\rm Hom}}
\def\hra{\hookrightarrow}

\def\i{\iota}

\def\Id{{\rm Id}}

\def\Im{{\rm Im}}

\def\Ker{{\rm Ker}}
\def\l{\lambda}
\def\la{\langle}

\def\L{\Lambda}
\def\lra{\longrightarrow}
\def\msk{\medskip}
\def\n{\noindent}

\def\o{\omega}
\def\op{\oplus}
\def\ot{\otimes}
\def\ov{\overline}
\def\O{\Omega}
\def\pf{\noindent{\bf Proof.}\hskip 2mm}

\def\PP{{\bf P}}

\def\qed{\hfill{\bf q.e.d.}}
\def\QQ{{\bf Q}}
\def\ra{\rangle}

\def\ri{\right}
\def\rk{{\rm rk}}

\def\RR{{\bf R}}
\def\s{\sigma}
\def\Si{\Sigma}

\def\ss{\subset}

\def\t{\theta}

\def\tm{\times}

\def\Tr{{\rm Tr}}

\def\ub{\underbrace}
\def\ul{\underline}
\def\vf{\varphi}

\def\wt{\widetilde}

\def\ZZ{{\bf Z}}
\magnification=1100
\centerline{\bf The weight-two Hodge  structure}

\centerline{\bf of moduli spaces of sheaves on a K3 surface.}
\bsk
\centerline{\bf Kieran G.~O'Grady}

\centerline{\bf October 1 1995}
\msk

\midinsert
\narrower\narrower\narrower
\n
{\bf Abstract.}
\hskip 2mm
We prove that the weight-two Hodge structure of moduli spaces of
torsion-free sheaves on a $K3$ surface  is as described by Mukai (the
rank is arbitrary but we assume the first Chern class is primitive). We prove
the moduli space is an irreducible symplectic variety (by Mukai's work  it
was known to be symplectic). By work of Beauville,  this implies that its
$H^2$ has a canonical integral non-degenerate quadratic form; Mukai's
recepee for $H^2$ includes a description of  Beauville's quadratic form. As
an application we compute higher-rank Donaldson polynomials of $K3$
surfaces.
\endinsert
\bsk

Recently Jun Li~[Li1,Li2]  determined the  stable rational Hodge
structure on the $n$-th cohomology, for $n\le 2$,  of moduli spaces of
rank-two torsion-free sheaves on an arbitrary projective surface (stable
means: for large enough dimension of the moduli space). It seems worthwile
to study in greater detail the integral Hodge structure of these moduli
spaces: in this paper we consider the case of a $K3$ surface.
Mukai~[M1,M2,M3] studied  extensively moduli of sheaves on  a
$K3$: in~[M2] he determined the  Hodge structure of
two-dimensional moduli spaces, and in~[M3, Th.~(5.15)] there is a beautiful
description of the  weight-two Hodge structure when the rank is at
most two, and the dimension of the moduli space is greater than two. In
this paper we  prove that Mukai's description of $H^2$ is valid  in arbitrary
rank. (But notice that we assume the first Chern class is primitive.)
\msk

\n
{\bf Statement of the results.}
\msk
\n
Let $S$ be a projective  $K3$ surface. The {\it Mukai lattice}~[M2, $\S$ 2]
consists of $H^*(S;\ZZ)$ endowed with the symmetric bilinear form
$$\,\,\langle\a,\b\rangle:=-\int_S\a^*\wedge\b\,,\eqno(1)$$
where for $\a=\a^0+\a^1+\a^2\in H^{*}(S)$, with
$\a^i\in H^{2i}(S)$, we let $\a^*:=\a^0-\a^1+\a^2$. Setting
$$\,\,F^0H^*(S):=H^*(S), \quad
F^1H^*(S):=H^0(S)\op F^1H^2(S)\op  H^4(S), \quad
F^2H^*(S):=F^2H^2(S)\,,$$
$(H^*(S;\ZZ),F^{\bu})$ becomes a
weight-two Hodge structure, with the additional datum of Mukai's quadratic
form defined above. If  $F$ is a sheaf on $S$, following Mukai we
set
$$\,\,v(F):=ch(F)(1+\o)=\rk(F)+c_1(F)+\left(\chi(F)-r\right)\o\,,$$
where $\o\in H^4(S;\ZZ)$ is the fundamental class. Notice that, since the
intersection form of $S$ is even, $v(F)\in H^*(S;\ZZ)$, and of course
$v(F)\in F^1H^*(S)$. An element of $F^1\cap H^*(S;\ZZ)$ will be
called a {\it Mukai vector}. Let $H$ be an ample divisor on $S$, and $v$ be a
Mukai vector:  we let $\cM_v(S,H)$  be the moduli space of
Gieseker-Maruyama $H$-semistable torsion-free sheaves $F$ on $S$ such
that $v(F)=v$. (We abbreviate to $\cM_v(H)$ whenever the surface $S$ is
fixed.) Thus $\cM_v(H)$ is a projective scheme. Now assume   $H$
is  {\it $v$-stabilizing}, i.e.~that all sheaves parametrized by $\cM_v(H)$
are $H$-slope-stable (see Proposition~(II.1)). In this case   the moduli
space is smooth~[M1, Th.~(0.3)],  of dimension equal to
$$\,\,d(v):=2+\la v,v\ra\,.$$
Furthermore there exists   a quasi-tautological family of  sheaves~[M2,
Th.~(A.5)], i.e.~a sheaf $\cF$ on $S\tm\cM_v(H)$, flat over
$\cM_v(H)$, such that if  $[F]\in\cM_v(H)$ represents the isomorphism
class of $F$, then
$$\cF|_{S\tm[F]}\cong F^{\op\s}$$
for some $\s>0$. We can assume $\s$ is independent of $[F]$, and we will
denote it by $\s(\cF)$.  Letting $\pi\cl S\tm \cM_v\to S$ and  $\rho\cl S\tm
\cM_v\to\cM_v$ be the  projections, we set
$$\,\,\t_{\cF}(\a):=
{1\over\s(\cF)}\rho_*\left[ch(\cF)^*(1+\pi^*\o)\pi^*\a\right]_3\,,$$
where $\a\in H^*(S;\CC)$, and $[\cdot]_3$ is the  component belonging to
$H^6(S\tm\cM_v(H))$.  If $\cE$ is another quasi-tautological sheaf then
there exist vector bundles $\xi$, $\eta$ on $\cM_v(H)$ such that
$\cE\ot\rho^*\xi\cong\cF\ot\rho^*\eta$~[M2, Th.~(A.5)]. From
this it follows easily that the restriction of $\t_{\cF}$ to $v^{\bot}\ss
H^*(S;\CC)$ is independent of the choice of a quasi-tautological family:
thus we get a canonical {\it Mukai map}~[M3, (5.14)]:
$$\,\,\t_v\cl v^{\bot}\to H^2(\cM_v(H);\CC)\,.$$
We will prove the following

\proclaim (2) Main Theorem.
Keeping notation as above, assume that $v^1$, the component of $v$
belonging to $H^2(S;\ZZ)$, is  primitive. Let $H$ be a $v$-stabilizing ample
divisor on $S$. Assume that the expected dimension of $\cM_v(H)$ is
greater than $2$, i.e.~that $\langle v,v\rangle>0$.
Then:
\msk
\item{1.} The moduli space $\cM_v(H)$ is an
irreducible symplectic variety, i.e.~simply connected with a symplectic
form spanning $H^{2,0}$, deformation equivalent to a symplectic
projective birational model of $T^{[n]}$,
where $T$ is a projective $K3$ surface,  $n:=d(v)/2$, and  $T^{[n]}$ is the
Hilbert scheme parametrizing  length-$n$ subschemes of $T$.  In particular,
by Beauville~[B, Th.~5], $H^2(\cM_v(H);\ZZ)$ has a canonical integral
quadratic form.
\item{2.} The map $\t_v$ is an isomorphism of  integral Hodge
structures, and an isometry if $v^{\bot}$ is provided with  the
restriction of Mukai's form, and $H^2(\cM_v(H)$ is provided Beauville's
quadratic form.
\msk

\n
{\bf Comments.}
\item{$\bu$}
Mukai~[M2] proved that if the  dimension of $\cM_v(H)$ is two,
i.e.~$\langle v,v\rangle=0$, then  the map $\t_v$ induces an
isomorphism of Hodge structures between $v^{\bot}/\CC v$ and
$H^2(\cM_v(H);\CC)$ (the only hypothesis is that semistability is
equivalent to stability).  Our proof of the main theorem can be easily
adapted to this case (with the hypothesis that $v^1$ is primitive). However
in a certain sense Mukai's proof is preferable: it is conceptual and it
"explains" the definition of  $\t_v$ and of the Mukai lattice. A conceptual
proof of Theorem~(2)  would be very interesting: it might indicate how to
extend the result to arbitrary surfaces.
\item{$\bu$}
We have some restrictions on the choice of $v$: it is natural to expect that
the theorem holds whenever semistability coincides with stability.
\item{$\bu$}
The key observation in our proof of the main theorem is the following: if $S$
is an elliptic $K3$ and $v^1$ has degree one on the elliptic fibers, then the
moduli space is birational to $S^{[n]}$.  It is interesting to notice that one
can   use the main theorem to show that in general the moduli space
is $\ul{\rm not}$ birational to any $T^{[n]}$~[M3, p.167].
\msk

\n
{\bf Plan of the paper.}
\hskip 2mm
First we prove the main theorem when $S$ is an elliptic $K3$ with a
section, $v^1$ has degree one on the fibers of the elliptic fibration, and
the ample divisor $H$ is suitable in the sense of Friedman~[F1]. In this
case the moduli space is birational to $S^{[n]}$: one possible route
to this result would be to proceed as in~[F2, Th.~(3.14)],
i.e.~to construct every stable vector-bundle as an elementary
modification of a fixed rigid bundle. Since we also want to prove Item~(2)
of the main theorem, we will proceed differently. We will prove that
$\cM_v(H)$ is birational to a moduli space of rank one less; iterating  we
get down to rank one, and since the moduli space of rank-one torsion-free
sheaves on $S$ (with $c_2=n$) is isomorphic to $S^{[n]}$, we conclude that
$\cM_v(H)$ is birational to $S^{[n]}$. This procedure allows us to verify
Item~(2) by induction on the rank, the case of rank one being trivial; all of
this is done in Section~I. In Section~II we prove the theorem in general: the
idea is to deform $S$ to an elliptic $K3$ (this part  is similar to an
argument in~[GoHu]).  In Section~III we give an application of the main
theorem: we compute higher-rank Donaldson polynomials. The last
section is devoted  to the proof of some technical results on
polarizations.
\msk

\n
{\it Acknowledgments.}
\hskip 2mm
The present work owes a lot to the papers of Mukai~[M1, M2, M3].
\bsk

\n
{\bf  I. The case of an elliptic $K3$ surface.}
\msk
\n
In order to state the main result of this section we need some
preliminaries on the choice of a polarization.  For a surface $S$, let
$A(S)\ss H^{1,1}(S;\RR)$ be the {\it ample cone}, i.e.~the real convex cone
spanned by   Chern classes of ample divisors.

\proclaim  Definition.
Let $k>0$. A {\it $k$-wall} of $A(S)$
consists of the intersection $L^{\bot}\cap A(S)$, where $L$ is a divisor on
$S$ such that
$$\,\,-k\le L^2<0\,.\eqno({\rm I}.0.1)$$
An open {\it $k$-chamber} $\cC\ss A(S)$ is a connected component of the
complement of the union of all $k$-walls.  An ample divisor on $S$ is {\it
$k$-generic} if it does not belong to any $k$-wall.

Let $S$ be a $K3$, and $v\in H^*(S;\ZZ)$ be a Mukai vector. We will often
consider $|v|$-walls,where
$$\,\,|v|:={(v^0)^2\over 4}\la v,v\ra+{(v^0)^4\over 2}\,.\eqno({\rm I}.0.2)$$
Now we specialize to the case of an elliptic $K3$,
i.e.~a $K3$ surface $S$ together with a linear pencil $|C|$,   where $C\ss S$
is an elliptic curve. (Thus $|C|$ defines a morphism to $\PP^1$ with generic
fiber an elliptic curve.)

\proclaim  Definition.
Let $(S,|C|)$ be a an elliptic $K3$, and $k$ be a positive integer. An ample
divisor $H$ on $S$ is {\it $k$-suitable} if it is $k$-generic and if $C$ is in
the closure of the unique open chamber containing $H$.

Notice that $H$ is k-generic (suitable) if and only if any ample divisor in
$\QQ_{+}H$ is  k-generic (suitable); hence it makes sense to consider
$k$-suitable polarizations.  One can show in general that $k$-suitable
polarizations exist for any $k$. We will only need the following special
case.

\proclaim (I.0.3) Lemma.
Let $S$ be an elliptic $K3$ surface with a section $\Si$ of the elliptic
fibration, and such that $Pic(S)=\ZZ[\Si]\op\ZZ[C]$, where $C$ is an elliptic
fiber.  Let  $H\sim(a\Si+bC)$ be an ample divisor such that
$$\,\,{b\over a}\ge k+1\,.$$
Then $H$ is $k$-suitable.

\pf
Proving that $H$ is $k$-suitable is equivalent to proving that
$sign(L\cdot H)=sign(L\cdot C)$
for all divisors $L\in {\rm Pic}(S)$  satisfying~(I.0.1). (Notice that
$(L\cdot C)\not=0$ for any such $L$.) Clearly it suffices to test only the $L$
such that $(L\cdot C)>0$. Let $L\sim(x\Si+yC)$. Then $L\cdot C>0$ is
equivalent to $x>0$, and~(I.0.1)  reads
$$\,\,0<2x(x-y)\le k\,.\eqno(*)$$
Using our hypothesis and the positivity of $a,x$ we get
$$\,\,L\cdot H=(b-2a)x+ay=a\cdot\left[\left({b\over a}-1\ri)x-(x-y)\ri]
\ge a\cdot\left(kx-(x-y)\ri)\,.$$
Since~($*$) holds, and since $x$, $(x-y)$ are positive integers, we have
$kx\ge 2(x-y)$. By the above inequality we conclude that $(L\cdot H)>0$.
\qed
\msk

A divisor $D$ (or the first Chern class of a divisor) on an elliptic $K3$ is a
{\it numerical section} if $D\cdot C=1$, where $C$ is an  elliptic fiber. This
section is devoted to proving   the following result.

\proclaim (I.0.4) Theorem.
Let $S$ be an elliptic $K3$. Assume that the elliptic fibration has a
section, and  that $\rho(S)=2$.  Let $v\in H^*(S;\ZZ)$ be
a Mukai vector such that  $v^1$  is a numerical section and $d(v)>2$.
Let $H$ be a $|v|$-suitable ample divisor on $S$.  Then $\cM_v(H)$ is an
irreducible symplectic variety birational to $S^{[n]}$,
where $n:=d(v)/2$. Furthermore the map
$$\t_v\cl H^*(S;\CC)\to H^2(\cM_v(H);\CC)$$
is an isomorphism of  integral Hodge structures, and an isometry.
\msk

\n
{\bf I.1. Preliminaries.}
\msk
Let $V$, $T$ be schemes. A {\it family of sheaves on $V$ parametrized by
$T$} consists of a sheaf $\cF$ on $V\tm T$, flat over $T$. For $t\in T$ we
let $\cF_t:=\cF|_{V\tm\{t\}}$.

\proclaim (I.1.1).
{\rm We recall  some well-known facts about stable bundles on elliptic
curves~[A] which will be useful in the course of the section. Let $C$ be an
elliptic curve. Then up to isomorphism there exists a unique stable bundle
$V_r^L$ on $C$ of given rank $r$ and determinant $L$. If $L$ has degree one
the picture is particularly simple. By Hirzebruch-Riemann-Roch, Serre
duality and stability we have $h^0(V_r^L)=1$. A non-zero section of $V_r^L$
gives rise to a non-split exact sequence
$$\,\,0\to \cO_C\to V_r^L\to V_{r-1}^L\to 0\,.$$
Conversely, since $\dim\Ext^1(V_{r-1}^L,\cO_C)=1$, the bundle  $V_r^L$
can be constructed as  the unique non-split extension  above.}

\proclaim Cayley-Bacharach.
{\rm We collect  some well-known criteria
for local-freeness of extensions~[GrHa, 729-731] involving the
Cayley-Bacharach property. Let $S$  be a surface. Let $L$, $M$ be
line-bundles on $S$. Let  $X,Y$ be $0$-dimensional reduced subschemes of
$S$, and $I_X$, $I_Y$ be their ideal sheaves. We consider extensions
$$\,\,0\to I_X\ot L\to F\to I_Y\ot M\to 0\,.\eqno({\rm I}.1.2)$$}

\proclaim (I.1.3).
Let $P\in Y$. If there exists a section of $L^{-1}\ot M\ot K_S$ vanishing at
all points of $(Y-P)$ and non-zero  at $P$,
then~(I.1.2) is split in a neighborhood of $P$, and hence $F$ is singular at
$P$.

\n
Now  assume that the image of the evaluation map
$$e_Y\cl H^0(L^{-1}\ot M\ot K_S)\to H^0(L^{-1}\ot M\ot K_S|_Y)$$
has codimension one. If~(I.1.2) is non-split at one point of $Y$, the
following converse of~(I.1.3) holds.

\proclaim (I.1.4).
Suppose $\Im(e_Y)$ has codimension one, and assume  that
Extension~(I.1.2) is non-split in the neighborhood of one point (at least)
of $Y$.  Let $P\in Y$. If all the sections of $L^{-1}\ot M\ot K_S$ vanishing
at $(Y-P)$ vanish also at $P$,  then $F$ is locally-free at $P$.

\proclaim (I.1.5).
Suppose that $h^1(I_X\ot M^{-1}\ot L)=0$. Let hypotheses be as
in~(I.1.4), except that instead of assuming~(I.1.2) is non-split in
the neighborhood of one point of $Y$, we only assume~(I.1.2) is globally
non-split. Then the conclusion of~(I.1.4) hold.

\n
{\bf Proof of~(I.1.3)-(I.1.4)-(I.1.5).}
\hskip 2mm
To prove the three statements  consider the exact sequence
$$\displaylines{\quad H^1(Hom(I_Y\ot M, I_X\ot L))\to
\Ext^1\left(I_Y\ot M, I_X\ot L\ri)\brel g\over\to\hfill\cr
\hfill{}\brel g\over\to
H^0\left(Ext^1(I_Y\ot M, I_X\ot L)\ri)
\brel f\over\to H^2(Hom(I_Y\ot M, I_X\ot L))\,.\quad\cr}$$
Identifying the dual of the last term with $H^0(M\ot L^{-1}\ot K_S)$ via
Serre duality, and the dual of $H^0(Ext^1(I_Y\ot M,I_X\ot L))$ with
$H^0(M\ot L^{-1}\ot K_S)|_Y)$ via Grothendieck duality, the    transpose of
$f$ gets identified with $e_Y$. From this~(I.1.3) follows at once. To
prove~(I.1.4) let
$$e_{(Y-P)}\cl H^0(L^{-1}\ot M\ot K_S)\to
H^0(L^{-1}\ot M\ot K_S|_{(Y-P)})$$
be evaluation. The obvious map $\Im(e_Y)\to \Im(e_{(Y-P)})$ is
an injection by hypothesis, and since $\Im(e_Y)$ has codimension one we
conclude that $e_{(Y-P)}$ is surjective. Let $\tau\in H^0(Ext^1(I_Y\ot M,
I_X\ot L))$ be the image under $g$ of the extension class corresponding
to~(I.1.2); since $e_Y$ is the transpose of $f$, $\tau$ is annihilated by
$\Im(e_Y)$. If the extension is split at $P$, we actually have $\tau\in
H^0(Ext^1(I_{(Y-P)}\ot M,I_X\ot L))$, and since $e_{(Y-P)}$ is surjective we
conclude that $\tau=0$, contradicting the assumption that~(I.1.2) is
non-split at one point of $Y$. Finally~(I.1.5) follows
from~(I.1.4) because if  $h^1(I_X\ot M^{-1}\ot L)$ vanishes then $g$ is
an injection.
\qed
\msk

The following proposition states a remarkable property of suitable
polarizations: the validity of this result is the reason for introducing
the notion of suitability. The proof will be given in Section~IV.

\proclaim (I.1.6) Proposition.
Let $S$ be an elliptic $K3$, with $|C|$ the elliptic fibration. Let $v\in
H^*(S;\ZZ)$ be a Mukai vector such that $\la v^1,C\ra$ and $v^0$ are
coprime, and let $H$ be a $|v|$-suitable ample divisor on $S$.
\msk
\item{\rm 1.} If a torsion-free sheaf $F$ with $v(F)=v$ is
$H$-slope-semistable then it is $H$-slope-stable. In particular $H$ is
$v$-stabilizing.
\item{\rm 2.} If $[F]\in\cM_v(H)$  the restriction of $F$ to the generic
$C_t\in|C|$ is stable.
\item{\rm 3.} Conversely if $F$ is a torsion-free sheaf on $S$ with
$v(F)=v$,  such that the restriction of $F$ to the generic
$C_t\in|C|$ is stable, then $F$ is $H$-slope-stable.
\msk

\n
{\bf I.2. Outline of the section.}
\msk
\n
Let $S$ be a $K3$ surface.

\proclaim (I.2.1) Definition.
Two Mukai vectors $v,w\in H^*(S;\ZZ)$  are {\it equivalent} ($v\sim w$) if
there exists a line bundle $\xi$ such that $w=ch(\xi)\cdot v$.

Thus if $F$ is a sheaf on $S$ then $v(F\ot\xi)\sim v(F)$.
Assume $v\sim w$: since multiplication by $ch(\xi)$ is
an isometry of the Mukai lattice we have $|v|=|w|$ (see~(I.0.2)), hence  open
$|v|$-chambers coincide with open $|w|$-chambers. From now on we assume
$S$ is elliptic with elliptic pencil $|C|$: if $v\sim w$  a polarization is
$|v|$-suitable if and only if it is $|w|$-suitable. If $H$ is $|v|$-suitable
the
map
$$\cM_v(H)\ni [F]\mapsto [F\ot\xi]\in\cM_w(H)$$
is an isomorphism: in fact if $[F]\in\cM_v(H)$ then $F$ is slope-stable
 by Proposition~(I.1.6), and
slope-stability  is preserved by tensorization.  Furthermore, since
multiplication by $ch(\xi)$ is an isometry of the Mukai lattice,
Theorem~(I.0.4) holds for $\cM_v(H)$ if and only if it holds for $\cM_w(H)$.
Thus we are allowed to replace $v$ by any equivalent vector; we will use
this freedom to normalize the Mukai vector of Theorem~(I.0.4) as follows.
Let $[F]\in\cM_v(H)$; since $v^1$ is a numerical section
$$\,\,\chi\left(F\ot[kC]\right)=\chi(F)+k\,.$$
Set
$$\hb{$\chi(v):=\chi(F)$, where $[F]\in\cM_v(H)$,}$$
and let $w:=v\cdot ch\left(\left[\left(1-\chi(v)\ri) C\ri]\right)$. Then
$w^1$ is again a numerical section, and furthermore  $\chi(w)=1$. Therefore
it suffices to prove Theorem~(I.0.4) under the additional hypothesis  that
$\chi(v)=1$: in this case we say $v$ is $\ul{\it normalized}$.

The proof of Theorem~(I.0.4) goes roughly as follows. We can assume $v$ is
normalized; if $[F]\in\cM_v(H)$ then by stability $h^2(F)$  vanishes and
hence $h^0(F)\ge\chi(v)=1$. Suppose that $h^0(F)=1$,    then we have a
canonical sequence
$$\,\,0\to\cO_S\to F \to E\to 0\,.$$
Assume also that $E$ is torsion-free (this will be the case if $F$ is
locally-free): since $H$ is $|v|$-suitable, it will follow that $E$ is
$H$-slope-stable. Let $\cM_w(H)$ be the moduli space to which the
isomorphism class of $E$ belongs: by mapping $[F]$ to $[E]$ we get a
rational map $\vf$ from $\cM_v(H)$ to $\cM_w(H)$.  Since $\chi(E)=-1$, we
see that  $h^1(E')\ge 1$ for all $[E']\in \cM_w(H)$. It turns out that
$h^1(E')=1$ for the generic $[E']\in\cM_w(H)$. Since
$h^1(E')=\dim\Ext^1(E',\cO_S)$  this means that $\vf$
has degree one, and thus $\cM_v(H)$ is
birational to $\cM_w(H)$. Normalizing $w$ and repeating this argument one
gets down to rank one, i.e.~$S^{[n]}$. This is the method by which we will
define a birational map between $\cM_v(H)$ and $S^{[n]}$; the details are in
the next subsection. In order to prove that $\t_v$ is an isomorphism of
Hodge structures we will  construct a subset $\cU$ of $\cM_v(H)$, and a
tautological family $\cF$ of sheaves on $S$ parametrized by $\cU$. Since
the complement of $\cU$ has codimension two, one has that $H^2(\cU)\cong
H^2(\cM_v(H))$, and that $\t_v$ is determined by $\t_{\cF}$.   This is the
longest part of the proof, the difficulty being that we must perform
semistable reduction along certain divisors; it takes up Subsections~I.4-I.5.
Finally  we will be able to verify that $\t_v$ is
an isometry by carrying out a purely numerical computation: this is the
content of the last subsection.
\msk

\n
{\bf I.3. The moduli space is birational to $S^{[n]}$.}
\msk
\n
For the rest of this section we assume $S$, $v$, $H$ are  as in
Theorem~(I.0.4), except that we do not assume $d(v)>2$, only $d(v)\ge 0$.
Furthermore we suppose that $v$ is normalized. We will sistematically
omit $H$ from our notation. We let $|C|$ be the elliptic pencil, and $\Si$ be
its section. Set
$$\,\,r:=v^0=\hb{rank of $F$, for $[F]\in\cM_v$}\qquad
n=d(v)/2=\hb{half the dimension of $\cM_v$}\,.$$
Then
$$v= r+\Si+(n-r^2+r)C+(1-r)\o\,.\eqno({\rm I}.3.1)$$
In particular  $v$ is determined by $r$, $n$;  we will denote $\cM_v$  by
$\cM^{2n}_r$.  The main result of this subsection is the following.

\proclaim (I.3.2) Proposition.
If  non-empty   the moduli space $\cM_r^{2n}$ is an irreducible
symplectic variety birational to $S^{[n]}$.

\n
In the next subsection we will show that in fact $\cM_r^{2n}$ is always
non-empty.  Before proving Proposition~(I.3.2) we need to discuss certain
Brill-Noether  loci. Let $W_v^d\ss\cM_v$ be the  subset
parametrizing locally-free sheaves $F$ such that $h^0(F)=(d+1)$ (we
assume $d\ge 0$).

\proclaim (I.3.3) Proposition.
Keep notation and assumptions as above.  Then $W^d_v$ is of pure
dimension, and
$$\,\,\cod(W_v^d,\cM_v)=\max\{d+1-\chi(v),0\}\,.$$

\n
We will prove Proposition~(I.3.3) at the end of this subsection.
\msk

\n
{\bf Proof of Proposition~(I.3.2).}
\hskip 2mm
It suffices to prove that $\cM_r^{2n}$  is  birational to $S^{[n]}$: in fact
since $\cM_r^{2n}$ is  symplectic~[M1], and since $S^{[n]}$ is symplectic
irreducible, it will follow that also $\cM_r^{2n}$ is symplectically
irreducible (i.e.~simply connected with a symplectic form
spanning $H^{2,0}$). The proof is by induction on $r$. When $r=1$ we have
$\cM_1^{2n}\cong S^{[n]}$, hence there is nothing to prove. So let's
suppose $r>1$. First we notice that the (open) subset of $\cM_r^{2n}$
parametrizing locally-free sheaves is dense: as is well-known this follows
from the fact that for all $[F]\in\cM_r^{2n}$ we have $h^2(\ad F^{**})=0$.
By  Proposition~(I.3.3) the open  subset $W_r^0:=W_v^0$ is dense in
$\cM_r^{2n}$.  If $[F]\in W_r^0$ then $F$ fits into a unique exact sequence
$$\,\,0\to \cO_S\to F\to Q\to 0\,.$$
The quotient $Q$ is torsion-free by the following.

\proclaim (I.3.4) Claim.
Keep notation as above.
Let $[F]\in \cM_v$, with $F$ locally-free, and suppose that
$h^0(F)=k+1>0$. Then $F$ fits into a unique exact sequence
$$\,\,0\to\cO_S(kC)\brel\a\over\to F\to Q\to 0\,,\eqno({\rm I}.3.5)$$
where $Q$ is torsion-free.

\n
{\bf Proof of the claim.}
\hskip 2mm
Let $f\cl S\to\PP^1$ be the elliptic fibration. Since $F$ is torsion-free
$f_*F$ is also torsion-free, hence locally-free. By~(I.1.1) $f_*F$  is of rank
one, and thus it is a line-bundle; since $h^0(f_*F)=h^0(F)$ we must have
$f_*F=\cO_{\PP^1}(k)$. The natural map $f^*f_*F\to F$ gives rise
to~(I.3.5), which is clearly unique. To prove that $Q$ is torsion-free
we must show that the divisorial part of the zero-locus of $\a$, call it
$div(\a)$, is empty. Let $C_t$ be a generic elliptic fiber; by~(I.1.6)
the restriction of $F$ to  $C_t$ is stable, and hence $div(\a)\cap C_t=\es$.
 Since all the elliptic fibers are irreducible ($\rho(S)=2$), we conclude that
 $div(\a)$ is a union of elliptic fibers. But if $\a$ vanishes on
an elliptic fiber  then we get a (non-zero) map
$\cO_S\left((k+1)C\right)\to F$, contradicting the assumption
$h^0(F)=k+1$. We conclude that $div(\a)=\es$, and thus $Q$ is torsion-free.
\qed
\msk

We go back to the proof of Proposition~(I.3.2). Let
$$\,\,w:=v(Q)=v(F)-1-\o\,.$$
By~(I.1.6) and by~(I.1.1) the
restriction of $Q$ to a generic elliptic fiber is stable, and hence by~(I.1.6)
we conclude that $Q$ itself is stable.
Thus $[Q]\in\cM_w$. Since $w^1=v^1(F)$, it is a numerical section and
hence $\cM_w$ is one of the moduli spaces we are considering. (But we do
not normalize $w$ for the moment.) From $h^0(F)=1$ we get  $h^0(Q)=0$.
By Serre duality $H^2(Q)\cong \Hom(Q,\cO_S)^*$, and by
stability of $Q$ the last group is zero. Thus $\chi(Q)=-h^1(Q)$. Since
$\chi(Q)=\chi(F)-\chi(\cO_S)=-1$ we conclude that $h^1(Q)=1$. Hence
$[Q]\in\cA_w$, where $\cA_w\ss\cM_w$ is the open subset
$$\,\,\cA_w:=\{[Q']\in\cM_w|\ h^1(Q')=1\}=\{[Q']\in\cM_w|\ h^0(Q')=0\}\,.$$
To sum up: we have defined a map
$$\,\,\vf\cl W_r^0\to\cA_w\,.$$
As is easily checked this map is a morphism. Let $[Q']\in\Im\vf$; since by
Serre duality $H^1(E)\cong\Ext^1(E,\cO_S)^*$, the morphism $\vf$ is
injective.  We claim that $\Im\vf$ is an open non-empty subset of $\cA_w$.
Let $[F]\in W_r^0$ ($[F]$ exists because by
hypothesis $\cM_r^{2n}$ is non-empty) and set $[Q]=\vf\left([F]\ri)$. By
definition of $\cA_w$ we have  $h^1(Q')=1$ for all $[Q']\in\cA_w$,
and since $\dim\Ext^1(Q',\cO_S)\cong H^1(Q')^*$ there exists a unique
non-trivial extension, call it $F'$,  of $Q'$ by $\cO_S$. By openness of
stability the sheaf $F'$  is stable for $[Q']$ varying in an open (non-empty)
subset of $\cA_w$; this proves that  $\Im\vf$ is an open non-empty subset
of $\cA_w$.  Let's  show that $\cA_w$ is dense in $\cM_w$.
It suffices to show that
$$\,\,\dim\{[Q']\in\cM_w|\ h^0(Q')>0\}
<\dim\cM_w\,.$$
If $w^0=(r-1)>1$ the locus parametrizing locally-free sheaves is
dense in $\cM_w$, and hence the  inequality follows from
Proposition~(I.3.3). If $w^0=1$ one easily checks the  inequality by hand.
Since $\vf$ is injective, and since its image is dense in $\cM_w$, it
defines a birational map between $\cM_r^{2n}$ and $\cM_w$.
Normalizing $w$ one gets $\cM_w\cong\cM_{r-1}^{2n}$.
By the inductive hypothesis $\cM_{r-1}^{2n}$  is birational to $S^{[n]}$, and
hence so is  $\cM_r^{2n}$.
\qed
\msk

\n
{\bf Proof of Proposition~(I.3.3).}
\hskip 2mm
First suppose $\chi(v)\ge(d+1)$. Let $[F]\in\cM_v$.
By Serre duality and stability $H^2(F)\cong \Hom(F,\cO_S)^*=0$, and hence
$h^0(F)\ge\chi(v)$. Thus $W^d_v$ is an open (eventually empty) subset of
$\cM_v$, and hence the proposition holds in this case (the empty set has
any codimension). From
now on we  assume that $\chi(v)\le d$.  One can describe  $W_v^d$ as a
determinantal variety: this is possible by  standard methods because
$H^2(F)$ vanishes for $[F]\in\cM_v$. The dimension formula for
determinantal varieties gives that
$$\,\,\cod(W_v^d,\cM_v)\le
\max\{(d+1)(d+1-\chi(v)),0\}\,.\eqno({\rm I}.3.6)$$
We first deal with a special case.

\proclaim (I.3.7) Lemma.
Keeping notation as above, suppose that $d=0$ and $\chi(v)\le 0$.    Then
$W_v^0$ is smooth of pure dimension, and
$$\,\,\cod(W_v^0,\cM_v)=1-\chi(v)\,.$$

\pf
Let $[F]\in W_v^0$. The non-zero section of $F$ gives an exact sequence
$$\,\,0\to\cO_S\to F\to Q\to 0\,.\eqno({\rm I}.3.8)$$
By Claim~(I.3.4) $Q$ is torsion-free. There is an exact sequence~[O,
Prop.~(1.17)]
$$\,\,0\to T_{[F]}W_v\to\Ext^1(F,F)\brel\b\over\to\Ext^1(\cO_S,Q)\,.$$
To compute $\rk\b$ consider the exact sequence
$$\,\,\Ext^1(F,F)\brel\g\over\to \Ext^1(F,Q)\to\Ext^2(F,\cO_S)
\to\Ext^2(F,F) \to\Ext^2(F,Q)\,.\eqno(*)$$
By Serre duality $\Ext^2(F,Q)\cong\Hom(Q,F)^*$. This last group is
zero: in fact $Q$ is stable because its restriction to a generic elliptic
fiber is stable, hence $\Hom(Q,Q)\cong\CC$, and since~(I.3.8) is not split we
conclude that $\Hom(Q,F)=0$.  Furthermore
$$\,\,\Ext^2(F,F)\cong\CC\qquad
\Ext^2(F,\cO_S)\cong\Hom(\cO_S,F)^*\cong\CC\,,$$
and hence the map $\g$ in~$(*)$ is surjective. From the exact
sequence
$$\,\,\Ext^1(F,Q)\brel\eta\over\to\Ext^1(\cO_S,Q)\to
\Ext^2(Q,Q)\to\Ext^2(F,Q)=0\,.$$
we get $dim\Im\eta=h^1(Q)-1$. Since $\b=\eta\circ\g$ we
conclude that
$$\,\,\cod(T_{[F]}W_v^0,T_{[F]}\cM_v)=h^1(Q)-1\,.$$
Since $h^0(Q)=h^2(Q)=0$, we have $h^1(Q)=-\chi(Q)$. From~(I.3.8) we get
$\chi(Q)=\chi(F)-2$, and hence
$$\,\,\cod(T_{[F]}W_v,T_{[F]}\cM_v)=1-\chi(F)=1-\chi(v)\,.$$
The lemma follows at once from this equality together with
Inequality~(I.3.6).
\qed
\msk

Now we prove Proposition~(I.3.3) in general. Let $[F]\in\cM_v$ with $F$
locally-free.  By Claim~(I.3.4) we have $h^0(F(-dC))=1$ if and only if  $[F]\in
W_v^d$. Hence, if $w:=v\cdot ch([-dC])$, tensorization by $[-dC]$ defines
an isomorphism between $W_v^d$ and $W_w^0$, and also of course between
$\cM_v$ and $\cM_w$. Thus
$$\,\,\cod(W_v^d,\cM_v)=\cod(W_w^0,\cM_w)
=1-\chi(w)=1-\chi(v)+d\,,$$
where the second equality holds by Lemma~(I.3.7).
\qed

\msk
\n
{\bf I.4. A large open subset of $\cM_r^{2n}$, for $r\le 2$.}
\msk
\n
Given $r\ge 1$ and $n\ge 0$ we will construct an open non-empty subset
$\cU_r^{2n}\ss\cM_r^{2n}$ and a tautological family of sheaves on $S$
parametrized by $\cU_r^{2n}$. Together with
Proposition~(I.3.2) this will establish that $\cM_r^{2n}$ is an irreducible
symplectic variety birational to  $S^{[n]}$. Furthermore since the
complement of $\cU_r^{2n}$ has codimension at least two the map $\t_v\cl
H^*(S)\to H^2(\cM_r^{2n})$ will be completely determined by a tautological
family on $S\tm\cU_r^{2n}$.   The construction of $\cU_r^{2n}$ and the
relative tautological family is by induction on $r$: the idea is to imitate
the picture for stable vector bundles on an elliptic curve (see~(I.1.1)). First
we will deal with the cases $r=1,2$. For simplicity's sake we fix $n$ and we
often omit it from our notation:  $\cM_r^{2n}$ will be denoted by
$\cM_r$, etc. The moduli space $\cM_1$ will be tacitly identified with
$S^{[n]}$.  Let $f\cl S\to\PP^1$ be the elliptic fibration.

\proclaim (I.4.1) Definition.
$\cU_1\ss S^{[n]}$ is the set consisting of $[Z]$ such that:
\msk
\item{1.}
$\# Z_{red}\ge (n-1)$, and $Z\cap \{\hb{critical points of $f$}\}=\es$,
\item{2.}
$h^0(I_Z\ot [(n-2)C])=0$, and if $h^0(I_Z\ot [(n-1)C])>0$ then $Z$ is
reduced,
\item{3.}
if $Z\cap \Si\not=\es$ then $Z$ is reduced, $Z\cap \Si$ consists of a
single  point,   and $h^0(I_Z\ot [(n-1)C])=0$.

\proclaim (I.4.2) Remark.
The complement of $\cU_1$ in $S^{[n]}$ has codimension two.

Our choice of a tautological family on $S\tm\cU_1$ (recall that we
identify $\cM_1$ with $S^{[n]}$) is
$$\,\,\cF^1:=I_{\cZ}\ot\pi^*[\Si+nC]\,,$$
where $\cZ\ss S\tm \cU_1$ is the tautological subscheme. Notice that if
$x\in\cU_1$ then $v(\cF^1_x)$ is the normalized Mukai vector $v$ with
$v^0=1$, $\la v^1,C\ra=1$, $d(v)=2n$, i.e.~the vector $v$ of~(I.3.1) with
$r=1$. Thus $\cF^1$ is indeed a tautological family of sheaves on
$S\tm\cU_1$. The reason for restricting to the subset $\cU_1$ of the whole
rank-one moduli space $S^{[n]}$ will become apparent when we deal with
higher-rank moduli spaces.  Now let's move to the case of rank two. We
will   construct a family of stable rank-two sheaves parametrized by
$\cU_1$, this family will define a classifying morphism $\cU_1\to\cM_2$,
and  $\cU_2$ will be defined as the image of this morphism.  The first step
is to construct a family of  extensions on $S$  parametrized by $\cU_1$. If
$[Z]\in\cU_1$ then by Items~(2)-(3) of~(I.4.1) together with Serre duality
$$\eqalignno{\,\,\dim\Ext^1(I_Z\ot[\Si+(n-2)C],\cO_S) &
=h^1(I_Z\ot[\Si+(n-2)C])=1\,,\cr
\Ext^i(I_Z\ot[\Si+(n-2)C],\cO_S) &
=0 \hb{ for $i=0,2$.} & ({\rm I}.4.3)\cr}$$
Hence if  $\rho\cl S\tm\cU_1\to\cU_1$ is projection,
$$\xi_2:=Ext^1_{\rho}(\cF^1\ot\pi^*[-2C],\cO_{S\tm\cU_1})$$
is a line-bundle on $\cU_1$.  The exact sequence
$$\displaylines{\quad 0=H^1(Ext^0_{\rho}(\cdot,\cdot))\to
\Ext^1\left(\cF^1\ot\pi^*[-2C]\ot\rho^*\xi_2,\cO_{S\tm\cU_1}\ri)\to
\hfill\cr
\hfill{}\to
H^0(Ext^1_{\rho}(\xi_2^{-1}\ot\cF^1\ot\pi^*[-2C],\cO_{S\tm\cU_1}))
\to H^2(Ext^0_{\rho}(\cdot,\cdot))=0\quad\cr}$$
shows that there is a unique non-trivial tautological extension
$$\,\,0\to\cO_{S\tm\cU_1}\to \cE^2\to
\cF^1\ot\pi^*[-2C]\ot\rho^*\xi_2\to 0\,.\eqno({\rm I}.4.4)$$
Since $\cO_{S\tm\cU_1}$ and $\cF^1$ are families of
torsion-free sheaves parametrized by $\cU_1$, so is $\cE^2$. If
$[Z]\in\cU_1$ then $\cE^2_{[Z]}$ is the unique non-split extension
$$\,\,0\to\cO_S\to \cE^2_{[Z]}\to
I_Z\ot[\Si+(n-2)C]\to 0\,.\eqno({\rm I}.4.5)$$
A computation shows that

\proclaim (I.4.6).
Keeping notation as above,  $v(\cE^2_{[Z]})$ is normalized and
$$\,\,v^0=2, \qquad\la v^1,C\ra=1,
\qquad 2+\la v(\cE^2_{[Z]}),v(\cE^2_{[Z]})\ra=2n\,.$$
Thus $v(\cE^2_{[Z]})$ is the vector  $v$ of~(I.3.1) with $r=2$.

If $\cE^2_{[Z]}$ is stable then by the above computation the isomorphism
class of $\cE^2_{[Z]}$  is represented by a point of $\cM_2$.  Before
analyzing  stability we prove the following proposition.

\proclaim (I.4.7) Proposition.
Let $Q$ be a torsion-free stable sheaf on $S$, with
$\la c_1(Q),C\ra=1$. Suppose that for $\ul{\rm all}$ $t\in\PP^1$ we have
$$\,\,\Hom(Q|_{C_t},\cO_{C_t})=0\,.\eqno({\rm I}.4.8)$$
Let
$$0\to[kC]\to F\to Q\to 0\eqno({\rm I}.4.9)$$
be a non-split extension.
\msk
\item{\rm 1.}
The sheaf $F$ is slope-stable.
\item{\rm 2.}
Letting
$$\,\,s:=\dim\Ext^1(Q,[kC])-1\,,$$
there are exactly $s$ elliptic fibers $C_u$
(counted with appropriate multiplicities) such that the restriction
of~(I.4.9) to $C_u$ splits.

\pf
By~(I.4.8) the sheaf  $Hom_f(Q,[kC])$ is zero.
{}From the exact sequence
$$0\to H^1\left(Hom_f\left(Q,[kC]\ri)\ri)\to \Ext^1\left(Q,[kC]\ri)\to
H^0\left(Ext^1_f(Q,[kC]\ri)\to H^2\left(Hom_f\left(Q,[kC]\ri)\ri)$$
we conclude that
$$\,\,\Ext^1(Q,[kC])\cong
H^0\left(Ext_f^1(Q,[kC])\right)\,.\eqno({\rm I}.4.10)$$
Since $Q$ is torsion-free it is $f$-flat,  hence the Euler
characteristic $\chi(Q|_{C_t},[kC]|_{C_t})$
is independent of the elliptic fiber $C_t$, and thus~(I.4.8)
implies that  $dim\Ext^1(Q|_{C_t},[kC]|_{C_t})$ is constant. Now let $C_t$ be
a generic elliptic fiber. Since $Q$ is stable its restriction to $C_t$ is
stable by~(I.1.6). By~(I.1.1) we get $dim\Ext^1(Q|_{C_t},[kC]|_{C_t})=1$, and
hence $\l:=Ext^1_f(Q,[kC])$ is  a line-bundle. By
Equality~(I.4.10) the extension class of~(I.4.9) corresponds to a non-zero
section of this line-bundle, and the elliptic fibers $C_u$ for
which the restriction of~(I.4.9) to $C_u$ is trivial are in one-to-one
correspondence with the zeroes of this section. This implies
Item~(1), by ~(I.1.1)-(I.1.6). It also implies  Item~(2) if we
notice that
$$\,\,\dim\Ext^1(Q,[kC])=h^0(\l)=\deg\l+1\,.$$
\qed

\proclaim (I.4.11) Corollary.
Let $[Z]\in\cU_1$. If $Z\cap\Si=\es$ and
$h^0\left(I_Z\ot[(n-1)C]\ri)=0$ then $\cE^2_{[Z]}$ is stable.

\pf
The sheaf $\cE^2_{[Z]}$ fits into the the non-split extension~(I.4.5).
Under our assumptions
$$\Hom(I_Z\ot[\Si+(n-2)C]|_{C_t},\cO_{C_t})=0\eqno(*)$$
for all elliptic fibers $C_t$. By Proposition~(I.4.7) we conclude that
$\cE^2_{[Z]}$ is stable.
\qed

\msk
Since for $[Z]\in\cU_1$ we have $\ell(Z)=n$, the above corollary shows
that $\cE^2_{[Z]}$ is stable for the generic $[Z]$. We let $L$ be the divisor
on $\cU_1$ consisting of points $[Z]$ such that
$$\,\,h^0\left(I_Z\ot[(n-1)C]\ri)>0\,,$$
i.e.~such that $Z$ contains two points lying on the same elliptic fiber.
If $D\ss S$ is an effective reduced divisor  we let $D_1$ be the
reduced divisor on $S^{[n]}$ defined by
$$\,\,D_1:=\{[Z]\in S^{[n]}|\ Z\cap D\not=\es\}\,.$$
Since the complement of $\cU_1$ in $S^{[n]}$ has codimension
two~(I.4.2) we can identify the group of divisors on $\cU_1$ with the
group of divisors in $S^{[n]}$; we will  use the same symbol for
corresponding divisors on  $\cU_1$ and $S^{[n]}$. The corollary above states
that $\cE^2_{[Z]}$ is stable if $[Z]\notin L\cup\Si_1$. We will show that
if $[Z]\in L\cup\Si_1$ then  $\cE^2_{[Z]}$ is indeed unstable.
\msk

\n
{\bf Semistable reduction along $L$.}
\hskip 2mm
Notice that if $n<2$ then $L=\es$, hence throughout this
subsubsection we assume that $n\ge 2$. Let $[Z]\in L$; since $[Z]\in L$, and
since by Item~(2) of~(I.4.1) the scheme $Z$ is reduced,  we can write
uniquely $Z=Z_0\cup W$, where  $Z_0$ consists of two points belonging to
the same  elliptic fiber $C_0$. To simplify notation we
set $E:=\cE^2_{[Z]}$. We recall that $E$ is the unique non-trivial extension:
$$\,\,0\to\cO_S\to E\to I_Z\ot[\Si+(n-2)C]\to 0\,.\eqno({\rm I}.4.12)$$
The sheaf $E$ is unstable for the following reason.  From
$$\dim\Hom(I_Z\ot[\Si+(n-2)C]|_{C_t},\cO_{C_t})=
\cases{0 & if $t\not=0$, and \cr
1 & if $t=0$,\cr}$$
it follows that
$\Ext^1_f(I_Z\ot[\Si+(n-2)C],\cO_S)\cong\cO_{\PP^1}(-1)\op\CC_0$,
where $\CC_0$ is the skyscraper sheaf at $0$.
Now notice that, setting $k=0$ and
$Q:=I_Z\ot[\Si+(n-2)C]$, Equality~(I.4.10) holds, because for its validity
it is sufficient that~(I.4.8) is satisfied for the generic elliptic
fiber. Hence we conclude that the restriction of the  non-split
extension~(I.4.12) to an elliptic fiber $C_t$ is split for $t\not=0$, and
non-split for $t=0$. (Strange as it may sound.) Hence the restriction of $E$
to the generic elliptic fiber is unstable, and by Proposition~(I.1.6) $E$
itself
must be unstable. It is also clear that a destabilizing sequence must have
some relation with the fiber $C_0$. Let's exhibit a destabilizing sequence.
Consider the natural exact sequence
$$\,\,0\to I_W\ot[\Si+(n-3)C]\brel\s\over\to
I_Z\ot[\Si+(n-2)C]\to \i_*\cO_{C_0}(-P)\to 0\,,\eqno({\rm I}.4.13)$$
where  $\i\cl C_0\hra S$ is the inclusion, and $P\in C_0$ is   the point such
that $(\Si\cap C_0)+P$ is linearly equivalent to $Z_0$ (as divisors on
$C_0$).  We claim that $\s$ lifts to a map $\wt{\s}\cl
I_W\ot[\Si+(n-3)C]\to E$. Consider the $\Ext$-sequence associated
to~(I.4.13):
$$\displaylines{\quad 0=\Hom(I_W\ot[\Si+(n-3)C],\cO_S)\to
\Ext^1(\i_*\cO_{C_0}(-P),\cO_S)\to\hfill\cr
\hfill{}\to\Ext^1(I_Z\ot[\Si+(n-2)C],\cO_S)
\brel\g\over\to \Ext^1(I_W\ot[\Si+(n-3)C], \cO_S)\,.\quad\cr}$$
The obstruction to lifting $\a$ is given by $\g(e)$, where $e$ is the
extension class of~(I.4.12). A computation shows that the first
$\Ext^1$-group appearing above is one-dimensional. Since the second
$\Ext^1$-group is also one-dimensional  we see that $\g=0$, and
hence $\s$ lifts to a map $\wt{\s}\cl I_W\ot[\Si+(n-3)C]\to E$.
The lift is unique because $h^0([-\Si-(n-3)C])=0$. We claim that the
quotient $E/\Im\wt{\s}$ is torsion-free. First notice that since $\s$ is an
isomorphism outside $C_0$, the quotient is certainly torsion-free outside
$C_0$. By~(I.1.5) $E$ is   locally-free along $C_0$, and clearly also
$I_W\ot[\Si+(n-3)C]$ is   locally-free along $C_0$. Thus it suffices to
show that $\wt{\s}$ does not vanish at the generic point of $C_0$. Assume
the contrary: we would get a non-zero map $I_W\ot[\Si+(n-2)C]\to E$,
which is absurd. Thus $E/\Im\wt{\s}$ is torsion-free; since its rank is one
it is isomorphic to  $I_Y\ot M$, where $Y$ is some zero-dimensional
subscheme of $S$, and $M$ is a line-bundle.  Computing Chern classes one
gets that $M=[C]$. Restricting $\wt{\s}$ to $C_0$ one gets that $Y=P$.
Hence we have an exact sequence
$$\,\,0\to I_W\ot[\Si+(n-3)C]\brel\b\over\to
E\to I_P\ot[C]\to 0\,.\eqno({\ rm I}.4.14)$$
This is  the Harder-Narasimhan filtration of $E$; since we do not
need this statement we omit its (easy) proof. By uniqueness of~(I.4.14)  we
can globalize the construction to all of $\cE^2|_{S\tm L}$:  letting $\cP\ss
S\tm L$ and $\cW\ss S\tm L$ be the subschemes swept out by $P$ and $W$
as $[Z]$ varies in $L$, there exist line bundles $\eta_1$, $\eta_2$ on $L$
such that we have an exact sequence
$$\,\,0\to I_{\cW}\ot\pi^*[\Si+(n-3)C]\ot\rho^*\eta_1\to \cE^2|_{S\tm L}
\to I_{\cP}\ot\pi^*[C]\ot\rho^*\eta_2\to 0\,.$$
Letting $\i^L\cl S\tm L\hra S\tm \cU_1$ be the
inclusion, we define $\cG^2$ to be the elementary modification of $\cE^2$
associated to the above exact sequence, i.~e.~$\cG^2$ is the sheaf
on $S\tm\cU_1$ fitting into the exact sequence
$$\,\,0\to \cG^2\to\cE^2\to
\i^L_*\left(I_{\cP}\ot\pi^*[C]\ot\rho^*\eta_2\ri)\to
0\,.\eqno({\rm I}.4.15)$$
The sheaf $\cG^2$ is $\cU_1$-flat~[F2, Lemma~(A.3)], i.e.~we can view it as
a family of sheaves on $S$ parametrized by $\cU_1$.

\proclaim (I.4.16) Proposition.
The sheaf $\cG^2_{[Z]}$ is torsion-free for all  $[Z]\in\cU_1$, and
stable if  $[Z]\in(\cU_1-\Si_1)$.

If  $[Z]\in\left(\cU_1-L\ri)$ then $\cG^2_{[Z]}\cong\cE^2_{[Z]}$, hence
the sheaf $\cG^2_{[Z]}$  is torsion-free  for
$[Z]\in\left(\cU_1-L\ri)$, and, by Corollary~(I.4.11),  stable for
$[Z]\in\left(\cU_1-L-\Si_1\ri)$. We are left with showing that
$\cG^2_{[Z]}$ is torsion-free and stable for $[Z]\in L$. Assume $[Z]\in L$;
for simplicity's sake set $G:=\cG^2_{[Z]}$.   The sheaf $G$ is naturally an
extension: in fact restricting~(I.4.15) to  $S\tm\{[Z]\}$ we obtain the
exact sequence
$$\,\,0\to I_P\ot [C]\to G
\to I_W\ot [\Si+(n-3)C]\to 0\,.\eqno({\rm I}.4.17)$$
Thus $G$ is torsion-free. Stability will be proved in various steps.
 First we consider the case when $n>2$.

\proclaim (I.4.18) Claim.
Keeping notation as above, assume $n>2$,  i.e.~$W\not=\es$.  There exists
at least one point of $W$ at which $G$ is locally-free.

Let's show that the above claim implies  $G$ is stable. Consider
$G^{**}$: from~(I.4.17) we obtain the exact sequence
$$\,\,0\to [C]\to G^{**}\to I_{W_0}\ot[\Si+(n-3)C]\to
0\,,\eqno({\rm I}.4.19)$$
where $W_0\ss W$ is the subset of points at which $G$ is locally-free
(recall that by~ Item~(2) of~(I.4.1) the scheme $W$ is reduced).
By~(I.4.18) we know $W_0\not=\es$, and hence the extension class
of~(I.4.19) is non-trivial. Since $W_0\cap\Si=\es$ and $W_0$ intersects
every elliptic fiber in at most one point, we can apply
Proposition~(I.4.7) to the extension~(I.4.19) and we conclude that $G^{**}$
is slope stable. Thus also $G$ is slope-stable.
\msk

\n
{\bf Proof of Claim~(I.4.18).}
\hskip 2mm
As usual  let $E:=\cE^2_{[Z]}$. The Kodaira-Spencer map
$$T_{[Z]}S^{[n]}\to\Ext^1(E,E)$$
is an isomorphism by the following criterion.

\proclaim (I.4.20).
Let $Q$ be a stable torsion-free  sheaf on $S$, such that
$c_1(Q)\cdot C=1$. Assume that
$$\,\,h^0(Q)=0  \qquad  h^1(Q)=\dim\Ext^1(Q,\cO_S)=1\,.\eqno(\dag)$$
Let $F$ be the unique non-trivial extension
$$\,\,0\to\cO_S\to F\to Q\to 0\,.\eqno(*)$$
Then $F$ is simple, and furthermore the Kodaira-Spencer  map
$\Ext^1(Q,Q)\to\Ext^1(F,F)$ corresponding to varying $Q$ and deforming
(uniqely by~($\dag$)) the extension~($*$) is an isomorphism.

\pf
The proof consists in  diagram chasing. We leave the details to the
reader.
\qed

\msk
The following is the key ingredient in the proof of Claim~(I.4.18).

\proclaim (I.4.21).
There exists at least one point $Q\in W$ such that the restriction
$$\,\,\Ext^1(E,E)\brel\phi_Q\over\to Ext^1(E_Q,E_Q)\,,$$
is surjective. (Here $E_Q$ is the localization of $E$ at $Q$.)

\pf
The local-to-global exact sequence for $\Ext$ gives an exact sequence
$$\,\,0\to H^1\left(Hom(E,E)\right)\to \Ext^1(E,E)
\brel\Phi\over\lra \bigoplus_{Q\in W}Ext^1(E_Q,E_Q)\,.\eqno(*)$$
Let's prove that
$$\,\,h^1\left(Hom(E,E)\right)=3\,.\eqno(**)$$
Considering the natural exact sequence
$$0\to Hom(E,E)\to Hom(E^{**},E^{**})\to
\bigoplus_{Q\in W}\CC_Q\to 0$$
we get
$$\,\,h^1\left(Hom(E,E)\right)=\dim\Hom(E,E)-h^0(E^{*}\ot E^{**})
+(n-2)+h^1(E^{*}\ot E^{**})\,.\eqno(\dag)$$
Considering the destabilizing sequence for $E^{**}$,  one
proves that $h^0(E^{*}\ot E^{**})=(n-2)$.
 Applying Hirzebruch-Rieman-Roch and Serre duality to $E^{*}\ot E^{**}$ we
get $h^1(E^{*}\ot E^{**})=2$. Since by~(I.4.20) we have
$\dim\Hom(E,E)=1$, Equation~($\dag$) gives~($**$). Now let's prove~(I.4.21),
arguing by contradiction. For all $Q\in W$ we have $E_Q\cong\cO\op I_Q$,
where $I_Q$ is the maximal ideal at $Q$, and hence   $\dim
Ext^1(E_Q,E_Q)=3$. Assuming that for all $Q\in W$ the restriction map
$\phi_Q$ is not surjective we conclude from~($*$) and~($**$) that
$\dim\Ext^1(E,E)<2n$. This is absurd because  $\dim\Ext^1(E,E)= 2n$ (by
Hirzebruch-Riemann-Roch and simplicity of $E$).
\qed

\msk
Now we are ready to prove Claim~(I.4.18). Let $Q\in W$ be a point
such that~(I.4.21) holds at $Q$. By~(I.4.20) and~(I.4.21) the map
$$\,\,T_{[Z]}S^{[n]}\to Ext^1(E_Q,E_Q)\,,\eqno(**)$$
obtained by composing  Kodaira-Spencer with restriction to
$Ext^1(E_Q,E_Q)$, is surjective, i.e.~the map  from a
neighborhood of  $[Z]$ to the versal deformation space
of $E_Q\cong I_Q\op\cO$ induced by $\cE^2$  is a submersion.  This means
that there are local coordinates $x,y$ on $S$, centered at $Q$,  and  local
coordinates $u,v,t$  on $S^{[n]}$, centered at $[Z]$, such that $\cE^2$ is
locally generated by $\a,\b,\g$ with the single relation
$$\,\,(x-u)\a+(y-v)\b+t\g=0\,.$$
By~(I.1.3)-(I.1.5) the sheaf $\cE^2_{[Z']}$, for $[Z']$ near $[Z]$,  is
singular exactly when $[Z']\in L$, and then it is singular at each point of
$W$. Hence a local equation for $L$ is given by $\{t=0\}$. Changing local
generators of $\cE^2$ if necessary we can assume that $\a$, $\b$ generate,
for $t=0$, the destabilizing subsheaf  $I_W\ot[\Si+(n-3)]$.    Then the
elementary modification $\cG^2$ (see~(I.4.15)) is locally  generated by $\a$,
$\b$ and $t\g$. By the above relation these elements generate a
locally-free (rank-two) sheaf.
\qed

\msk
This finishes the proof of Proposition~(I.4.16) when $n>2$. Now we examine
the case $n=2$.  The destabilizing sequence~(I.4.14) becomes
$$\,\,0\to [\Si-C]\to E\to I_P\ot[C]\to 0\,,\eqno({\rm I}.4.22)$$
and $E$ is locally-free. Exact sequence~(I.4.17) becomes
$$\,\,0\to I_P\ot[C]\to G\to [\Si-C]\to 0\,.\eqno({\rm I}.4.23)$$
Let $C_t$ be a generic elliptic fiber: we will prove   $G|_{C_t}$ is
stable. This will establish that $G$ is stable by Proposition~(I.1.6). The
vector bundle  $G|_{C_t}$ is stable if and only if the restriction
of~(I.4.23) to $C_t$ is non-split. To determine the relevant extension
class we proceed as follows.  Consider  $\cE^2_{[Z']}|_{C_t}$ for
$[Z']\in\cU_1$: by~(I.4.22)  this is unstable for $[Z']\in L$, while by~(I.1.6)
it is stable (for $C_t$  generic) if $[Z']\in(\cU_1-L-\Si_1)$.  Furthermore,
restricting Exact sequence~(I.4.15) to $C_t\tm\cU_1$, we see that
$\cG^2|_{C_t\tm\cU_1}$ is obtained from $\cE^2|_{C_t\tm\cU_1}$ by
applying the first step of semistable reduction to the bundles parametrized
by $L$. Hence by~[O, (1.11)] the extension class of the restriction of~(I.4.23)
to $C_t$ is given by the following recipe. Let
$$\a\cl T_{[Z]}S^{[2]}\to H^1\left(\cO_{C_t}(-\Si)\ri)$$
be the composition of Kodaira-Spencer and  the natural map  $H^1(\ad E)\to
H^1\left(\cO_{C_t}(-\Si)\right)$ induced by restriction of Exact
sequence~(I.4.22). If  $x\in T_{[Z]}S^{[2]}$ is a tangent vector, the class
$\a(x)$ represents the obstruction to lifting to first order, in the direction
$x$, the restriction
 of~(I.4.22) to $C_t$. Thus $\a$ vanishes on $T_{[Z]}L$, and hence it
induces a map
$$\,\,\ov{\a}\cl T_{[Z]}S^{[2]}/T_{[Z]}L\to
H^1\left(\cO_{C_t}(-\Si)\right)\,.$$
The extension class of the restriction of~(I.4.23) to $C_t$ is given by a
generator of $\Im\ov{\a}$. By~(I.4.20) the Kodaira-Spencer
map is surjective, and thus it suffices to show that
$$\,\,H^1(\ad E)\to H^1\left(\cO_{C_t}(-\Si)\right)\eqno({\rm I}.4.24)$$
is non-zero. Consider
$$\,\,H^1(\ad E)\to H^1(\ad E|_{C_t})\brel\b\over\to H^2(\ad E\ot[-C])\,.$$
By Serre duality the transpose of $\b$ is identified with the restriction
map
$$\,\,^{t}\b\cl H^0\left(\ad E\ot[C]\ri)\to H^0(\ad E|_{C_t})\,.$$
Let's examine $H^0\left(\ad E\ot[C]\ri)$. Using~(I.4.22) one
verifies  that any $\vf\in \Hom(E,E\ot[C])$ decomposes as
$\vf=\vf^0+\vf^1$, where $\vf^0\in H^0([C])\ot\Id_E$,  and
$\vf^1$ is a constant multiple of the composition
$$\,\,E\to I_P\ot[C]\to[C]\to E\ot[C]\,,$$
where the last map is given by the non-zero section of $E$. Thus
 $H^0(\ad E\ot[C])$ is one-dimensional; let $\vf$ be a
generator. Then
$$\,\,\vf|_{C_t}= \left[\matrix{\eta_{11} & \eta_{12} \cr
0 & -\eta_{11}\cr}\right]\,,$$
where $\eta_{11}\in H^0(\cO_{C_t})$, and $\eta_{12}\in
H^0(\cO_{C_t}(\Si)$ (recall that
$E|_{C_t}\cong\cO_{C_t}(\Si)\op\cO_{C_t}$),  and if $C_t$ is generic then
$\eta_{11}\not=0$.  Now let $\tau\in H^1(\ad E|_{C_t})$. Then
$$\,\,\tau= \left[\matrix{\tau_{11} & 0 \cr
\tau_{21} & -\tau_{11}\cr}\right]\,,$$
where $\tau_{11}\in H^1(\cO_{C_t})$, and $\tau_{21}\in
H^1(\cO_{C_t}(-\Si))$. Since $\Im^{t}\b$ is genereted by $\vf|_{C_t}$,
Serre duality tells us that $\tau\in\Ker\b$ if and only if
$$\,\,0=\Tr\left[\matrix{\tau_{11} & 0 \cr
\tau_{21} & -\tau_{11}\cr}\right]\cdot
\left[\matrix{\eta_{11} & \eta_{12} \cr
0 & -\eta_{11}\cr}\right]=
2\tau_{11}\eta_{11}+\tau_{21}\eta_{12}\,.$$
Since $\eta_{11}\not=0$ we conclude that the generator $\tau$ of $\Ker\b$
has $\tau_{21}\not=0$. Let $\wt{\tau}\in H^1(\ad E)$ be a lift of $\tau$.
The image of $\wt{\tau}$ under the map~(I.4.24) is  equal to $\tau_{21}$,
which is non-zero. Hence $\Im\ov{\a}$ is non-zero, and  $G|_{C_t}$
is stable. This proves  Proposition~(I.4.16) in the case $n=2$, and
finishes the proof of the proposition.
\msk

\n
{\bf Semistable reduction along $\Si_1$.}
\hskip 2mm
Notice that $\Si_1$ is empty if $n=0$, hence we assume throughout this
subsubsection that $n\ge 1$. Let $[Z]\in\Si_1$. By Item~(3) of
Definition~(I.4.1) the intersection $Z\cap\Si$ consists of a single reduced
point $P$, and $Z=P\cup W$, where $W$ is reduced of length $(n-1)$. We
let $C_0$ be the elliptic fiber through $P$.    To simplify notation we set
$G:=\cG^2_{[Z]}$.  Since $G$ is the same as $\cE^2_{[Z]}$, it is presented as
the unique non-trivial extension
$$\,\,0\to\cO_S\to G\to I_Z\ot[\Si+(n-2)C]\to 0\,.\eqno({\rm I}.4.25)$$
The sheaf $\cG^2_{[Z]}$ is unstable; one argues exactly as when we
explained why $\cE^2_{[Z]}$ is unstable for $[Z]\in L$. We will not go
through the same argument again. Instead we directly exhibit a destabilizing
sequence for $G$.   Let
$$\,\,0\to I_W\ot[\Si+(n-3)C]\brel\s\over\to
I_Z\ot[\Si+(n-2)C]\to \i_*\cO_{C_0}\to 0\,,\eqno({\rm I}.4.26)$$
be the natural exact sequence, where   $\i\cl C_0\hra S$ is the inclusion.
We claim that $\s$ lifts to a map $\wt{\s}\cl I_W\ot[\Si+(n-3)C]\to G$.
Consider the $\Ext$-sequence associated to~(I.4.26):
$$\displaylines{\quad 0=\Hom(I_W\ot[\Si+(n-3)C],\cO_S)\to
\Ext^1(\i_*\cO_{C_0},\cO_S)\to\hfill\cr
\hfill{}\to\Ext^1(I_Z\ot[\Si+(n-2)C],\cO_S)
\brel\g\over\to \Ext^1(I_W\ot[\Si+(n-3)C], \cO_S)\,.\quad\cr}$$
The obstruction to lifting $\a$ is given by $\g(e)$, where $e$ is the
extension class of~(I.4.25). Since the first $\Ext^1$-group
in the above exact sequence has dimension zero, and
since also the second $\Ext^1$-group is one-dimensional, we see
that $\g=0$, and hence $\s$ lifts to a map $\wt{\s}\cl
I_W\ot[\Si+(n-3)C]\to G$. The lift is unique because $h^0([-\Si-(n-3)C])=0$.
The quotient $G/\Im\wt{\s}$ is torsion-free: the argument is the same
as the one given to show that $E/\Im\wt{\s}$ is torsion-free when we
performed semistable reduction along $L$, we omit to repeat it.
Since $E/\Im\wt{\s}$ is torsion-free of rank  one,
it is isomorphic to  $I_Y\ot M$, where $Y$ is some zero-dimensional
subscheme of $S$ and $M$ is a line-bundle.  Computing Chern classes one
gets that $M=[C]$ and $Y$ is empty. Hence we get
$$\,\,0\to I_W\ot[\Si+(n-3)C]\to G\to [C]\to 0\,.\eqno({\rm I}.4.27)$$
This sequence gives the Harder-Narasimhan
filtration of $G$.   By unicity  the above exact sequence
globalizes: thus there are line-bundles $\t_1$, $\t_2$ on $\Si_1$ such that
we have
$$\,\,0\to I_{\cW}\ot\pi^*[\Si+(n-3)C]\ot\rho^*\t_1\to
\cG^2|_{S\tm\Si_1}\to\pi^*[C]\ot\rho^*\t_2\to 0\,,$$
where $\cW$ is the scheme swept out by $W$ as $[Z]$ varies in $\Si_1$.
Let $\cF^2$ be the sheaf on $S\tm\cU_1$ defined by
$$\,\,0\to\cF^2\to\cG^2\to
\i_*^{\Si_1}\left(\pi^*[C]\ot\rho^*\t_2\right)\to 0\,,\eqno({\rm I}.4.28)$$
where $\i^{\Si_1}\cl \Si_1\hra\cU_1$ is the inclusion. The sheaf $\cF^2$
is $\cU_1$-flat ([F2, Lemma~(A.3)]), hence we can view it as a family of
sheaves on $S$ parametrized by $\cU_1$.

\proclaim (I.4.29) Proposition.
Keep notation as above. If $[Z]\in\cU_1$ then
$v\left(\cF^2_{[Z]}\right)$ is the  normalized  Mukai
vector $v$ with $v^0=2$, $\la v^1,C\ra=1$, and $d(v)=2n$, i.e.~it is equal to
the vector  of~(I.3.1) with $r=2$.  The sheaf $\cF^2$ is a family of
stable torsion-free sheaves on $S$, parametrized by $\cU_1$.

\pf
Since $\cF^2$ is $\cU_1$-flat, and since $\cU_1$ is connected, the vector
$v\left(\cF^2_{[Z]}\right)$ is independent of $[Z]$. If
$[Z]\in(\cU_1-L-\Si_1)$ then $\cF^2_{[Z]}=\cE^2_{[Z]}$, thus the first
statement follows from~(I.4.6).
For $[Z]\in(\cU_1-\Si_1)$ we have $\cF^2_{[Z]}=\cG^2_{[Z]}$, hence the
only thing left to prove is that $\cF^2_{[Z]}$ is stable and locally-free for
$[Z]\in\Si_1$.  Restricting~(I.4.28) to $S\tm[Z]$, for $[Z]\in\Si_1$, one
gets
$$\,\,0\to [C]\to\cF^2_{[Z]}\to I_W\ot[\Si+(n-3)C]\to
0\,.\eqno({\rm I}.4.30)$$ %
For each $t\in\PP^1$ one
has $\Hom\left(I_W\ot[\Si+(n-3)C]|_{C_t},[C]|_{C_t}\right)=0$, because of
Item~(3) of~(I.4.1). Thus by Proposition~(I.4.7) the sheaf $\cF^2_{[Z]}$ is
stable if the extension~(I.4.30) is non-split. The extension class of~(I.4.30)
is obtained as follows. Let
$$\a\cl T_{[Z]}S^{[n]}\to \Ext^1\left(I_W\ot[\Si+(n-3)C],[C]\right)$$
be the composition of Kodaira-Spencer for $\cG^2$ and of the  map
$$\b\cl \Ext^1\left(\cG^2_{[Z]},\cG^2_{[Z]}\right)\to
\Ext^1\left(I_W\ot[\Si+(n-3)C],[C]\right)$$
 associated to ~(I.4.27).
For $x\in T_{[Z]}S^{[n]}$ the class $\a(x)$ is the obstruction to
lifting~(I.4.27) to first order in the direction of $x$. Thus $\a$ vanishes
on $T_{[Z]}\Si_1$, and it induces a map
$$\,\,\ov{\a}\cl T_{[Z]}S^{[n]}/T_{[Z]}\Si_1\to
\Ext^1\left(I_W\ot[\Si+(n-3)C],[C]\right)\,.$$
The extension class of~($*$) is equal to the generator of $\Im\ov{\a}$~[O,
(1.11)]. By~(I.4.20) the Kodaira-Spencer map is surjective, hence it suffices
to check that the map $\b$ is non-zero. One checks that the lemma below
applies to the exact sequence~(I.4.27), and hence $\b$ does not vanish.
This proves Extension~(I.4.30) is non-split, and hence  $\cF^2_{[Z]}$ is
stable by Proposition~(I.4.7).
\qed

\proclaim  Lemma.
Let $A$, $B$ be stable torsion-free sheaves on $S$, with
$$\,\,\dim\Ext^1(A,B)\ge 2\qquad\Hom(A,B)=\Hom(B,A)=0\,.$$
Suppose we have a non-trivial extension
$$\,\,0\to A\to F\to B\to 0\,.$$
Then the map $\Phi\cl \Ext^1(F,F)\to \Ext^1(A,B)$ is non-zero.

\pf
The proof consists of diagram chasing. First one proves that $F$ is simple,
and thus by Serre duality $\Ext^2(F,F)\cong\CC$. Then one gets surjectivity
of $\Ext^1(F,F)\to \Ext^1(F,B)$. Finally one considers
$$\Ext^1(F,B)\brel\Psi\over\to\Ext^1(A,B)\to
\Ext^2(B,B)\to \Ext^2(F,B)\cong\Hom(B,F)^*=0\,.$$
Since $B$ is stable it is simple, and thus $\dim\Ext^2(B,B)=1$.
Since $\dim\Ext^1(A,B)\ge 2$ the map $\Psi$ is non-zero, and hence
$\Phi$ does not vanish.
\qed

\msk
Since $\cF^2$ is a family of stable sheaves on $S$ parametrized by
$\cU_1$, with normalized Mukai vector $v$, it defines a classifying
morphism $\vf\cl\cU_1\to\cM_2$.

\proclaim (I.4.31) Proposition-Definition.
Keep notation as above. Let $\cU_2:=\vf(\cU_1)$. Then $\cU_2$ is an open
subset of $\cM_2$, and its complement has codimension two. The map
$\vf\cl \cU_1\to\cU_2$ is an isomorphism, and hence we can view $\cF^2$
as a tautological family of sheaves on $S$ parametrized by $\cU_2$.

\pf
On the subset $\left(\cU_1-L-\Si_1\ri)$ the map $\vf$ is injective. In fact
if $[F]\in\vf\left(\cU_1-L-\Si_1\ri)$ then $h^0(F)=1$ by Items~(2)-(3)
of~(I.4.1), and the inverse $\vf^{-1}$ on $\vf\left(\cU_1-L-\Si_1\ri)$ is
obtained by associating to $[F]$ the (zero-dimensional) zero-locus of a
non-zero section of $F$. Thus $\vf$ has degree one onto its image. Let $\t$
be a symplectic form on $\cM_2$~[M1]. Since the complement of
$\cU_1$ in $S^{[n]}$ has codimension two,  $\vf^*\tau$ extends
to a two-form on all of $S^{[n]}$. This two-form is non-zero because $\vf$
has degree one, and hence it is a symplectic form. This implies that $\vf$
is an embedding;  since $\dim\cU_1=\dim\cM_2$ and since $\cM_2$,
$\cU_1$ are smooth we conclude that $\vf$ is an isomorphism onto
$\cU_2$, and that $\cU_2$ is open in $\cM_2$. Now let's prove that the
complement of $\cU_2$ has codimension at least two. By
Proposition~(I.3.2) we know $\cM_2$ is irreducible, hence $\vf$ defines a
birational map
$$\,\,\vf\cl S^{[n]}\cdots>\cM_2\,.$$
Assume that $(\cM_2-\cU_2)$ contains a divisor $D$. The map
$\vf^{-1}$ is a morphism on some (non-empty) open subset $D^0\ss D$
because $S^{[n]}$, $\cM_2$ are both irreducible symplectic. Applying
Zariski's main theorem to a resolution of indeterminacies of $\vf$ we
see that $\vf^{-1}(D^0)\ss\left(S^{[n]}-\cU_1\ri)$. Since $(S^{[n]}-\cU_1)$
has codimension two, we conclude that $\vf^{-1}$ contracts $D^0$.
This is absurd because the pull-back by $\vf^{-1}$ of a symplectic form on
$S^{[n]}$ is a symplectic form on $\cM_2$. Thus $(\cM_2-\cU_2)$ has
codimension at least two in $\cM_2$.
\qed
\msk

\n
{\bf I.5. A large open subset of $\cM_r$, for $r\ge 3$.}
\msk
We retain  notations and conventions of the previous subsection.
For each $r\ge 3$ we will construct a subscheme $\cU_r$ of $\cM_r$ and
a tautological family $\cF^r$ of  sheaves on $S$
parametrized by $\cU_r$.
In order to construct $\cF^r$  we need  the
following inductive hypotheses:

\proclaim (I.5.1).
$\cU_r\cong \cU_1$.

\proclaim (I.5.2).
Let $r\ge 2$.  If $x\in\cU_r$ then $h^0(\cF_x^r\ot[-2C])=0$.

\proclaim (I.5.3).
Let $r\ge 2$, and let $y\in\Si_1\ss\cU_r$ (this makes sense
by~(I.5.1)).   Then $\cF^r_y$ is locally-free. There exists exactly
one elliptic  fiber $C_u$ such that $\Hom(\cF_y^r|_{C_u},\cO_{C_u})\not=0$.
Furthermore  $\cF_y^r|_{C_u}\cong \cO_{C_u}\op V$,
where  $\Hom(V,\cO_{C_u})=0$.

\proclaim (I.5.4).
Let $r\ge 2$. If $x\in(\cU_r-\Si_1)$ then
$\Hom(\cF_x^r|_{C_t},\cO_{C_t})=0$ for $\ul{\rm all}$ elliptic  fibers
$C_t$.

Let $r\ge 2$ be given, and assume~(I.5.2)-(I.5.4) hold for this  $r$:
we will construct $\cF^{r+1}$.   Let
$x\in\cU_r$. Since the Mukai vector corresponding to $\cM_r$ is
normalized we have $\chi(\cF^r_x)=1$. Let $C_s$, $C_t$ be two elliptic
fibers; the exact sequence
$$0\to \cF_x^r\ot[-2C]\to\cF^r_x\to\cF^r_x|_{C_s\cup C_t}\to 0$$
shows that $\chi(\cF_x^r\ot[-2C])=-1$. By~(I.5.2)
$h^0(\cF_x^r\ot[-2C])=0$, and  by Serre duality together with stability of
$\cF^r_x$,
$$\,\,H^2(\cF^r_x\ot[-2C])\cong
\Hom(\cF^r_x\ot[-2C],\cO_S)^*=0\,.\eqno({\rm I}.5.5)$$
We conclude that $H^1(\cF_x^r\ot[-2C])$ is one-dimensional. By Serre
duality $\Ext^1(\cF^r_x\ot[-2C],\cO_S)$ is one-dimensional.   Thus if
$\rho\cl S\tm\cU_r\to\cU_r$ is projection,
$$\xi_{r+1}:=Ext^1_{\rho}(\cF^r_x\ot[-2C],\cO_{S\tm\cU_r})$$
is a line-bundle on $\cU_r$. Let  $\cE^{r+1}$ be the tautological extension
$$\,\,0\to\cO_{S\tm\cU_r}\to\cE^{r+1}\to
\cF^r\ot\pi^*[-2C]\ot\rho^*\xi_{r+1}\to 0\,.$$
Since $\cF^r$ is $\cU_r$-flat, so is $\cE^{r+1}$. If $x\in \cU_r$ then
$\cE^{r+1}_x$ is the unique non-trivial extension
$$\,\,0\to\cO_S\to\cE_x^{r+1}\to\cF_x^r\ot[-2C]\to 0\,.
\eqno({\rm I}.5.6)$$

\proclaim (I.5.7) Lemma.
Let $r\ge 2$ be given, and assume~(I.5.4) holds for this $r$.
Then $\cE^{r+1}$ is a family of torsion-free sheaves on $S$ parametrized by
$\cU_r$. If $x\in\cU_r$,  $v(\cE_x^{r+1})$ is the
normalized Mukai vector $v$ such that
$$\,\,v^0=(r+1) \qquad \la v^1,C\ra =1 \qquad d(v)=2n\,,\eqno({\rm I}.5.8)$$
i.e.~the vector corresponding to $\cM_r$.
For $x\in(\cU_r-\Si_1)$ the sheaf $\cE^{r+1}_x$ is stable, and furthermore
for all elliptic fibers $C_t$ one has
$\Hom(\cE^{r+1}_x|_{C_t},\cO_{C_t})=0$.

\pf
$\cU_r$-flatness follows immediately from $\cU_r$-flatness of $\cF^r$.
Since $\cF^r_x$ is torsion-free, so is $\cE^{r+1}$. To show
that $v(\cE^{r+1}_x)$ is normalized consider the long exact cohomology
sequence associated to~(I.5.6). Equation~(I.5.8) is
proved by a simple computation.
Stability follows  from Item~(1) of Proposition~(I.4.7). By Item~(2) of
Proposition~(I.4.7) the
restriction of~(I.5.6) to any elliptic fiber is  non-split, and hence
the result follows from Remark~(I.5.9) below.
\qed

\proclaim (I.5.9) Remark.
Let $C_t$ be an elliptic fiber, and $0\to \cO_{C_t}\to V\to W\to 0$ be a
non-trivial extension. If $\Hom(W,\cO_{C_t})=0$ then
$\Hom(V,\cO_{C_t})=0$.

Let's show that for $y\in\Si_1$ the sheaf
$\cE^{r+1}_y$ is not stable.   By~(I.5.3)
there exists a unique elliptic fiber $C_u$ such that
$\Hom(\cF^r_y\ot[-2C]|_{C_u},\cO_{C_u})\not=0$, this group is
one-dimensional, and if $\b$ is a generator, the corresponding map is
surjective.    We define $\cH^r_y$ to be the sheaf on $S$ fitting
into the exact sequence
$$\,\,0\to \cH^r_y\brel\a\over\to\cF^r_y\ot[-2C]
\brel\b\over\to\i^u_*\cO_{C_u}\to 0\,,\eqno({\rm I}.5.10)$$
where $\i^u\cl C_u\hra S$ is the inclusion.
At this point we  add to our list of inductive hypotheses one  last item.

\proclaim (I.5.11).
Let $r\ge 2$ be given. Let $y\in\Si_1$. Then $\cH_y^r$
is locally-free, and $\Hom(\cH_y^r|_{C_t},\cO_{C_t})=0$ for all elliptic
fibers $C_t$.

We claim that $\a$ lifts to $\cE^{r+1}_y$; this will give a destabilizing
subsheaf of $\cE^{r+1}_y$. Consider the exact sequence
$$\,\,\Hom(\cH^r_y,\cO_S)\to \Ext^1(\i^u_*\cO_{C_u},\cO_S)\to
\Ext^1(\cF^r_y\ot[-2C],\cO_S)\brel\g\over\to\Ext^1(\cH^r_y,\cO_S)\,.$$
Letting $e$ be the extension-class of~(I.5.6),  $\g(e)$ is the obstruction to
lifting $\a$ to $\cE^{r+1}_y$.  Since  $\cH_y^r|_{C_t}=\cF_y^r|_{C_t}$ for
$t\not=u$, the restriction of $\cH_y^r|_{C_t}$ to the generic elliptic fiber
is stable, and hence  the first group of the above sequence is zero. A
straightforward computation shows that the first $\Ext^1$-group
is  one-dimensional; since we know that also the second $\Ext^1$-group
has dimension one, $\g$ must
vanish. This proves that $\a$ lifts to a map  $\wt{\a}\cl
\cH^r_y\to\cE^{r+1}_y$.

\proclaim Claim.
Let $r\ge 2$ be given. Assume that (I.5.2), (I.5.3)
and~(I.5.11) hold for $r$. Then
$$\,\,\cE^{r+1}_y/\Im\wt{\a}\cong [C]\,.$$

\pf
The quotient $\cE^{r+1}_y/\Im\wt{\a}$ is clearly of rank one. Let's start by
showing that it is torsion-free. By~(I.5.11) the sheaf $\cH^r_y$  is
locally-free. By~(I.5.3) $\cF^r_y\ot[-2C]$ is locally-free, hence by Exact
sequence~(I.5.6) also $\cE^{r+1}_y$ is locally-free.  Thus, since outside of
$C_u$ the map $\wt{\a}$ is an injection,  $\cE^{r+1}_y/\Im\wt{\a}$ has
torsion if and only if $\wt{\a}$ drops rank along all of $C_u$. Let's assume
that this is indeed the case: since
$\Ker\left(\wt{\a}|_{C_u}\ri)\ss\Ker\left(\a|_{C_u}\ri)$, and since the
second kernel is a line-bundle, we conclude that the two kernels coincide.
We claim this implies that $\wt{\a}$ extends to a map
$\cF^r_y\ot[-2C]\to\cE^{r+1}_y$. In fact consider the subsheaf ${\rm
Sat}(\Im\wt{\a})\ss\cE^{r+1}_y$ associated to the presheaf
$$\,\,{\rm Sat}(\Im\wt{\a})_U:=\{\s\in\left(\cE^{r+1}_y\ri)_U|
\hb{ $\vf\s\in\Im\wt{\a}$ for some $0\not= \vf\in\cO_U$}\}\,.$$
${\rm Sat}(\Im\wt{\a})$ is locally-free, and since
$\Ker\left(\wt{\a}|_{C_u}\ri)$ is of rank one,  it fits into an exact
sequence
$$\,\,0\to \cH^r_y\brel\wt{\a}\over\to{\rm Sat}(\Im\wt{\a})
\to\i^u_*\cO_{C_u}\to 0\,.$$
Since $\Ker\left(\wt{\a}|_{C_u}\ri)=\Ker\left(\a|_{C_u}\ri)$, we have
${\rm Sat}(\Im\wt{\a})\cong\cF^r_y\ot[-2C]$. Thus
$\wt{\a}$ extends to a map  $\cF^r_y\ot[-2C]\to\cE^{r+1}_y$. This is
impossible because~(I.5.6) does not split. Hence $\cE^{r+1}_y/\Im\wt{\a}$
is torsion-free, and therefore isomorphic to $I_Y\ot\l$, for a line-bundle
$\l$ and a zero-dimensional subscheme $Y$ of $S$. Computing Chern classes
one gets that $\l\cong[C]$ and that $Y=\es$.
\qed
\msk

Thus we get a unique exact sequence:
$$\,\,0\to\cH_y^r\to\cE^{r+1}_y\to[C]\to 0\,.\eqno({\rm I}.5.12)$$
Since $\la c_1(\cH_y^r), C\ra=1$, the above exact sequence shows that the
restriction of $\cE^{r+1}_y$ to the generic fiber is unstable; by~(I.1.6) we
conclude that $\cE^{r+1}_y$ is unstable.  By unicity of the above
construction, we can globalize it to all of $S\tm\Si_1$. Hence there exists
a line-bundle $\eta_{r+1}$ on $S\tm\Si_1$ and a surjection
$\cE^{r+1}|_{S\tm\Si_1}\to\eta_{r+1}$, which restrict on $S\tm\{y\}$ to
the above destabilizing quotient for $\cE^{r+1}_y$. We define   $\cF^{r+1}$
as the elementary modification of $\cE^{r+1}$ fitting into the exact
sequence
$$\,\,0\to\cF^{r+1}\to\cE^{r+1}
\to\i_*^{\Si_1}\eta_{r+1}\to 0\,.\eqno({\rm I}.5.13)$$

\proclaim (I.5.14) Lemma.
Let $r$ be given, and assume~(I.5.2)-(I.5.11) hold. Then $\cF^{r+1}$ is a
family of  torsion-free stable sheaves on $S$. If $x\in\cU_r$ then
$v(\cF_x^{r+1})$ is equal to the normalized Mukai vector $v$
satisfying~(I.5.8). If $y\in\Si_1$ there exists exactly one elliptic fiber
$C_u$ such that $\Hom(\cF_y^{r+1}|_{C_u},\cO_{C_u})$ is non-zero.
Furthermore  $\cF_y^{r+1}|_{C_u}\cong\cO_{C_u}\op V$, where
$\Hom(V,\cO_{C_u})=0$.

\pf
Since $\cE^{r+1}$ is $\cU_r$-flat, the sheaf $\cF^{r+1}$ is flat by~[F2,
Lemma~(A.3)].
$$\hb{If $x\in\cU_r-\Si_1$ then $\cF^{r+1}_x=\cE^{r+1}_x$.}
\eqno({\rm I}.5.15)$$
If $y\in\Si_1$ restriction of~(I.5.13) to $S\tm\{y\}$ gives
$$\,\,0\to[C]\to\cF_y^{r+1}\to\cH_y^r\to 0\,.\eqno({\rm I}.5.16)$$
In particular $\cF^{r+1}_y$ is torsion-free, hence $\cF^{r+1}$ is a
family of torsion-free sheaves on $S$ parametrized by $\cU_r$. From this
it follows that there  is a two-step locally-free resolution of $\cF^{r+1}$,
and hence the Chern character of $\cF^{r+1}_x$ is a locally constant
function of $x\in\cU_r$. Thus $v(\cF_x^{r+1})$  is a locally constant
function of $x\in\cU_r$. Since, by~(I.5.1), $\cU_r$ is irreducible, this
function is constant; hence Equality~(I.5.15) and Lemma~(I.5.7) imply that
$v(\cF_x^{r+1})$ is as stated. Let's prove stability.  For
$x\in(\cU_r-\Si_1)$ this follows from~(I.5.15) together with
Lemma~(I.5.7). Let $y\in\Si_1$. Proceeding exactly as in the proof of
Proposition~(I.4.29) one shows that Extension~(I.5.16) is non-trivial, and
hence stability of  $\cF_y^{r+1}$ follows from~(I.4.7) together
with~(I.5.11). What is left to  prove are  the last two statements. By~(I.5.11)
we can apply Proposition~(I.4.7) to Extension~(I.5.16): since
$\dim\Ext^1(\cH_y^r,[C])=2$ we get that there exists a unique elliptic fiber
$C_u$ such that the restriction of~(I.5.16) to $C_u$ is trivial. By
Remark~(I.5.9) together with~(I.5.11) we conclude that
$\Hom(\cF^{r+1}_y|_{C_t},\cO_{C_t})\not=0$ if and only if $t=u$.
That $\cF^{r+1}_y|_{C_u}$ splits as stated follows immediately
from~(I.5.11).
\qed
\msk

Now we can define $\cU_{r+1}$.  By Lemma~(I.5.14) the sheaf
$\cF^{r+1}$  defines a  classifying
morphism $\vf\cl\cU_r\to\cM_{r+1}$.

\proclaim (I.5.17) Proposition-Definition.
Keep notation as above. Let $\cU_{r+1}:=\vf(\cU_r)$. Then $\cU_{r+1}$ is an
open subset of $\cM_{r+1}$, and its complement has codimension two.  The
morphism $\vf\cl \cU_r\to\cU_{r+1}$ is an
isomorphism onto its image. Hence we can view  $\cF^{r+1}$ as a
tautological family of sheaves on $S$ parametrized $\cU_{r+1}$.

\pf
The proof goes exactly as the proof of Proposition~(I.4.31).
\qed
\msk

What is left to do is to prove inductively that
hypotheses~(I.5.1)-(I.5.11) hold for all $r\ge 2$. Item~(I.5.1) is
obvious.
\msk

\n
{\bf Proof of~(I.5.2).}
\hskip 2mm
To verify~(I.5.2) for $r=2$ one considers separately the three cases
$x\in(\cU_2-L-\Si_1)$, $x\in L$, and $x\in\Si_1$. Set $x=[Z]$. In the first
case $\cF^2_{[Z]}=\cE^2_{[Z]}$; tensoring Exact
sequence~(I.4.5) by $[-2C]$, and using Definition~(I.4.1), one
gets the desired result. In the second case $\cF^2_{[Z]}=\cG^2_{[Z]}$, and
one gets the result arguing similarly, with Sequence~(I.4.17)
replacing~(I.4.5).  In the third case one proceeds in the same way,
with~(I.4.30) replacing~(I.4.5).  Now let's prove the inductive step. If
$x\in(\cU_{r+1}-\Si_1)$ then $\cF^{r+1}_x=\cE^{r+1}_x$ (of course we
identify $\cU_{r+1}$ with $\cU_r$); tensoring~(I.5.6) by $\cO_S(-2C)$ and
applying the inductive hypothesis one gets the result. Let $y\in\Si_1$.
Tensoring~(I.5.16) by $\cO_S(-2C)$ we see that it suffices to show that
$h^0(\cH^r_y\ot[-2C])=0$. By~(I.5.10) $\cH^r_y$ is a subsheaf of
$\cF^r_y\ot[-2C]$ hence the result follows from the inductive hypothesis.
\msk

\n
{\bf Proof of~(I.5.3).}
\hskip 2mm
First we verify that~(I.5.3) holds for $r=2$. We start by showing that
$\cF^2_y$ is locally-free for all $y\in\Si_1$. The proof is by contradiction.
Assume $\cF^2_{[Z_0]}$ is singular for some $[Z_0]\in\Si_1$. Let's prove
that this implies $\cF^2_{[Z]}$ is singular for all $[Z]\in\Si_1$. Since
$\Si_1\cap L=\es$, there exists $[Z_0]\in(\cU_2-L)$ such that
$\cF^2_{[Z_0]}$ is singular. The locus of $[Z]\in(\cU_2-L)$ such that
$\cF^2_{[Z]}$ is singular has codimension at most one, because $\cF^2$ has
rank two.  If $[Z]\in\left(\cU_2-L-\Si_1\ri)$ then $\cF^2_{[Z]}=\cE^2_{[Z]}$;
since $\cE^2_{[Z]}$ fits into the non-trivial extension~(I.4.5) it follows
from~(I.1.5) together with Definition~(I.4.1) that $\cF^2_{[Z]}$ is
locally-free. Since $\Si_1$ is irreducible we conclude that $\cF^2_{[Z]}$ is
singular for $\ul{\rm all}$ $[Z]\in\Si_1$. Let $[Z]\in\Si_1$, and set
$Z=P\cap W$, where $P:=Z\cap\Si$.  Since $\cF^2_{[Z]}$ fits into exact
sequence~(I.4.30) its singular points belong to $W$. On the other hand we
claim that if $n>2$ (we will treat the case $n=2$ separately) there exists
at least one point of $W$ at which $\cF^2_{[Z]}$ is locally-free. In fact
consider the exact sequence
$$\displaylines{\quad H^1\left(Hom(I_W\ot[\Si+(n-3)C],[C]\ri)\to
\Ext^1\left(I_W\ot[\Si+(n-3)C],[C]\ri)\brel\tau\over\to\hfill\cr
\hfill{}\brel\tau\over\to
H^0\left(Ext^1(I_W\ot[\Si+(n-3)C],[C])\ri)
\to H^2\left(Hom(I_W\ot[\Si+(n-3)C],[C]\ri)\,.\quad({\rm I}.5.18)\cr}$$
Since $n>2$ the first term is zero, and the claim follows immediately. Now
fix $P_0\in\Si$, and consider
$$\,\,\O_0:=\{[W]\in\left(S-P_0\ri)^{[n-1]}|\ [P_0\cup W]\in\Si_1\}\,.$$
Each $[W]\in\O_0$ is the disjoint union of the subset consisting of
the points at which $\cF^2_{[W\cup P_0]}$ is locally-free and of its
complement. Since for each $[W]\in\O_0$  these subsets are
both non-empty, we conclude that the incidence locus
$$\wt{\O}_0:=\{\left([W],Q\ri)\in \O_0\tm S|\ Q\in W\}$$
has two components. This is absurd, and hence $\cF^2_{[Z]}$ is locally-free
for all $[Z]\in\Si_1$, if $n>2$. Now let's suppose that $n=2$. Assume
$\cF^2_{[Z]}$ is singular. Consider Exact sequence~(I.5.18): since $W$
consists of one point, and since
$$\,\,\dim H^0\left(Ext^1(I_W\ot[\Si+(n-3)C],[C])\ri)=1\,,$$
the image under $\tau$ of the extension class
of~(I.4.30) must be zero, i.e.~the extension class lives in the first
cohomology group of~(I.5.18). Since this group is one-dimensional, we
conclude that the isomorphism class of $\cF^2_{[Z]}$ only depends on $W$,
and not on $P$. This is absurd because $\cU_2$ is a subset of the moduli
space $\cM_2$. We have finished proving that  $\cF^2_{[Z]}$ is locally-free
for all $[Z]\in\Si_1$.   To prove the second and third statements
of~(I.5.3) one applies  Item~(2) of Proposition~(I.4.7) to the non-split
exact sequence~(I.4.30): this gives that there exists a unique elliptic fiber
$C_u$ such that the restriction of~(I.4.30) to $C_u$ splits. By
Remark~(I.5.9) and Definition~(I.4.1) one concludes that
$\Hom(\cF^2_{[Z]}|_{C_t},\cO_{C_t})\not=0$ if and only if $t=u$, and that
the splitting at $t=u$ is as claimed. It remains to verify the inductive step:
this was proved in Lemma~(I.5.14).
\msk

\n
{\bf Proof of~(I.5.4).}
\hskip 2mm
Let's prove~(I.5.4) for $r=2$. First
consider $x\in(\cU_2-\Si_1-L)$. Set $x=[Z]$. Then $\cF^2_{[Z]}=\cE^2_{[Z]}$.
By Item~(2) of Proposition~(I.4.7) the restriction of~(I.4.5) to any
elliptic fiber is non-split, and hence Statement~(I.5.4) follows from
Remark~(I.5.9). We are left with proving that~(I.5.4) holds for $x\in L$ (and
$r=2$).

\proclaim Claim.
Assume there exists $x\in L$ and an elliptic fiber $C_u$ such that
$\Hom(\cF_x^2|_{C_u},\cO_{C_u})\not=0$. Then $\cF^2_x|_{C_u}$ is
singular.

\n
{\bf Proof of the claim.}
\hskip 2mm
Since $x\in L$, we have $\cF^2_x=\cG^2_x$.
To simplify notation set  $G:=\cG_x^2$. Then we have
$$\,\,0\to I_P\ot[C]\to G\to I_W\ot[\Si+(n-3)C]\to 0\,,$$
with notation as in~(I.4.17). By~(I.4.18) the sheaf $G$
 is locally-free at one at least of  the points of $W$. By~(I.1.4)
and Definition~(I.4.1) the sheaf  $G$ is locally-free at all
points of $W$, and hence we get
$$\,\,0\to [C]\to G^{**}\to I_W\ot[\Si+(n-3)C]\to 0\,,\eqno(\dag)$$
Since $\Ext^1(I_W\ot[\Si+(n-3)C],[C])$ is one-dimensional, the restriction
of~($\dag$) to any elliptic fiber $C_t$ is non-split by Item~(2) of
Proposition~(I.4.7). Thus by
Remark~(I.5.9) $\Hom(G^{**}|_{C_t},\cO_{C_t})=0$ for all elliptic fibers
$C_t$. Since $G|_{(S-P)}\cong G^{**}|_{(S-P)}$ we conclude that $P\in C_u$,
i.e.~$G|_{C_u}$ is singular.
\qed
\msk

Let's go back to the proof that~(I.5.4) holds for $r=2$ and $x\in L$. Assume
that~(I.5.4) does not hold for some $x\in L$;  we will arrive at a
contradiction. Let $G$, $C_u$  be as above. Let
$\Def^0(G|_{C_u})$ be the subspace of the deformation space
$\Def(G|_{C_u})$ parametrizing deformations "fixing the determinant".  The
restriction map $r_u\cl \Def(G)\to \Def^0(G|_{C_u})$  is surjective. In fact
it suffices to check that $h^2(\ad F^{**}\ot [-C])=0$; by Serre duality this is
equivalent to $h^0(\ad F^{**}\ot[C])=0$, and this follows from stabilty of
$F^{**}|_{C_t}$ for a generic elliptic fiber $C_t$. Letting
$V\ss\Def^0(F|_{C_u})$ be the subspace parametrizing deformations $A$
such that $\Hom(A,\cO_{C_u})\not=0$, one verifies easily that the generic
point of $V$ parametrizes a locally-free sheaf. Since the restriction map
$r_u$ is surjective we conclude that there exists $x'\in\cU_2$ arbitrarily
near to $x$ such that  $\Hom(\cF_{x'}^2|_{C_u},\cO_{C_u})\not=0$ and
$\cF^2_{x'}|_{C_u}$ is locally-free. By the previous claim $x'\notin L$. Since
$\Si_1\cap L=\es$ we can assume $x'\in(\cU_2-L-\Si_1)$. This is absurd:
by Item~(I.5.3), which we have just proved,
$\Hom(\cF_{x'}^2|_{C_u},\cO_{C_u})\not=0$ for such an $x'$. Now  let's prove
the inductive step. Since $x\in(\cU_{r+1}-\Si_1)$ we have
$\cF^{r+1}_x=\cE^{r+1}_x$ (we are identifying $\cU_{r+1}$ with $\cU_r$).
By Remark~(I.5.9) and  Inductive hypothesis~(I.5.4) it suffices to show
that the restriction of~(I.5.6) to any elliptic fiber is non-split. This
follows at once from Item~(2) of Proposition~(I.4.7).
\msk

\n
{\bf Proof of~(I.5.11).}
\hskip 2mm
By~(I.5.3), which we have just proved, $\cF^r_y$ is locally-free for
$r\ge 2$. Exact sequence~(I.5.10) proves that also $\cH^r_y$ is locally-free
for $r\ge 2$. Now let's prove that
$$\hb{$\Hom(\cH^r_y|_{C_t},\cO_{C_t})=0$
for all elliptic fibers $C_t$.}\eqno({\rm I}.5.19)$$
For this it is convenient to notice that
$\cH^r_y$ has been defined also when $r=1$. In fact Exact sequence~(I.5.10)
in the case $r=1$ becomes  Exact sequence~(I.4.26). Thus
$$\,\,\cH^1_{[Z]}=I_W\ot[\Si+(n-3)C]\,.$$
Exact sequences~(I.5.12)-(I.5.16) for $r=1$  reduce to~(I.4.27)
and~(I.4.30) respectively. We will prove inductively that~(I.5.19) holds for
all $r\ge 1$. It holds for $r=1$ because by Definition~(I.4.1)
 $W$ intersects every elliptic fiber in at most one reduced point, not
belonging to $\Si$. Now assume~(I.5.19) holds for $r\ge 1$.    Replacing $r$
by $r+1$ in~(I.5.10), and  tensorizing the exact sequence  by $[2C]$,  we get
$$\,\,0\to\cH_x^{r+1}\ot[2C]\to\cF^{r+1}_x\to\i^u_*\cO_{C_u}\to 0\,.$$
Thus local sections of
$\cH_x^{r+1}\ot[2C]$ consist of local sections of $\cF^{r+1}_x$ whose
projection to $[C]|_{C_u}\cong\cO_{C_u}$ is zero, where the projection is
defined by the splitting of~(I.5.16) along $C_u$. Let $\s$ be a non-zero
section  of $[C]$ which vanishes on $C_u$: by~(I.5.16) this gives a
section of $\cF_x^{r+1}$, which we still denote by  $\s$.  Since $\s$
vanishes on $C_u$, it belongs to  $H^0\left(\cH_x^{r+1}\ot[2C]\ri)$; notice
that as a section of this locally-free sheaf it is nowhere zero along $C_u$.
Considering~(I.5.16) we see that the resulting
exact sequence is:
$$\,\,0\to\cO_S\brel\s\over\to\cH_x^{r+1}[2C]
\to\cH_x^r\to 0\,.\eqno({\rm I}.5.20)$$
The restriction of the above extension to a generic elliptic fiber
coincides with the restriction of~(I.5.16), and hence it is non-trivial. By
the inductive hypothesis $\Hom(\cH_x^r|_{C_t},\cO_{C_t})=0$ for all
elliptic fibers. Thus we can apply Item~(2) of Proposition~(I.4.7) to
Extension~(I.5.20). A computation  gives $\dim\Ext^1(\cH_x^r,\cO_S)=1$,
and hence by~(I.4.7) the restriction of the
above extension to any elliptic fiber is non-trivial. By~(I.5.9) and
the inductive hypothesis we conclude that~(I.5.11) holds with
$r$ replaced by $(r+1)$.
\msk

\n
This ends the proof that the tautological families $\cF^r$ are defined for
all $r$.
\msk

\n
{\bf I.6. Proof of Theorem~(I.0.4).}
\msk
As explained in Subsection~(I.2) we can assume $v$ is
normalized, that is $\cM_v(H)=\cM_r^{2n}$ for $r=v^0$ and $n=d(v)/2$. By
Subsections~I.4-I.5 we know $\cM_r^{2n}$ is non-empty,  thus by
Proposition~(I.3.2) it is birational to $S^{[n]}$. We are left with proving
that $\t_v$ is an isomorphism of integral  Hodge structures, and that it
preserves the quadratic forms. As in the previous subsections we think $n$
is fixed, and  we omit it from our notation whenever this shouldn't cause
confusion. We let
$$v_r:=r+c_1(\Si+(n-r^2+r)C)+(1-r)\o$$
be the Mukai vector corresponding to $\cM_r^{2n}$, and we set
$\t_r:=\t_{v_r}$. A {\it lattice} $(\L,q)$ consists of  a free $\ZZ$-module
$\L$ provided with an integral quadratic form $q$.  A homomorphism
between lattices $(\L_1,q_1)$, $(\L_2,q_2)$ consists of a
homomorphism of modules $f\cl\L_1\to\L_2$ such that
$q_1(\a)=q_2(f(\a))$ for all $\a\in\L_1$.

\proclaim (I.6.1) Proposition.
Keep notation as above and assume $n\ge 1$. The map $\t_r$ is
 integral, i.e.
$$\,\,\t_r\left(v_r^{\bot}\cap H^*(S;\ZZ)\ri)\ss H^2(\cM_r;\ZZ)\,.$$
This map is a homomorphism of lattices, where $v_r^{\bot}\cap H^*(S;\ZZ)$
is provided with Mukai's quadratic form, given by~(1), and
$H^2(\cM_r;\ZZ)$ is provided with Beauville's canonical quadratic form
$B_{\cM_r}$~[B, Th.~(5)].

\n
{\bf Proof of Proposition~(I.0.4) assuming~(I.6.1).}
\hskip 2mm
Our first step is to identify the quadratic form
$B_{\cM_r^{2n}}$. For this we need the lemma below. This result  is
well-known to experts~[M3, Prop.~(5.8)]; we give a proof for the
reader's convenience.

\proclaim (I.6.2) Proposition.
Let $X,Y$ be birational irreducible symplectic projective varieties. Let
$f\cl X\cdots>Y$ be a birational map.
There exist open subsets $i\cl U\hra X$, $j\cl V\hra Y$, whose
complements have codimension at least two, such that $f$ restricted to
$U$ is regular and it defines an isomorphism $f|_U\cl U\brel\cong\over\to
V$. Let $f^{\sharp}\cl H^2(X;\ZZ)\brel\cong\over\to H^2(Y;\ZZ)$ be the
isomorphism obtained as the composition
$$\,\,H^2(X;\ZZ)\brel i^*\over{\brel\sim\over\rightarrow}
H^2(U;\ZZ)\brel\sim\over\to  H^2(V;\ZZ)
\brel j^*\over{\brel\sim\over\leftarrow} H^2(Y;\ZZ)\,,$$
which is independent of $U$, $V$. Then for all $\a\in H^2(X)$
$$\,\,B_X(\a)=B_Y(f^{\sharp}\a)\,.\eqno({\rm I}.6.3)$$

\pf
Let $X\leftarrow Z\rightarrow Y$ be a  resolution of
indeterminacies of $f$, with $Z$ smooth. Let's
show that the exceptional divisors of $g$ are the same as the exceptional
divisors of $h$. If $\tau_X\in H^0(\O_X^2)$ then there exists $\tau_Y\in
H^0(\O_Y^2)$ such that
$$\,\,g^*\tau_X=h^*\tau_Y\,.\eqno(*)$$
Now assume $\tau_X\not=0$. Let $\dim X=\dim Y=2n$. If $E$ is an
exceptional divisor of $g$ then $\wedge^n(g^*\tau_X)$ vanishes on $E$,
hence by~($*$) so does $\wedge^n(h^*\tau_Y)$. Since $\tau_Y$ is
non-degenerate this implies that $E$ is an exceptional divisor of $h$. Thus
every exceptional divisor of $g$ is an exceptional divisor of $h$. Reversing
the r\^oles of $X$ and $Y$ we conclude that the exceptional divisors of
$g$ and $h$ are the same. Clearly the first
statement holds if we set $U:=g(\cup_i E_i)$, $V:=h(\cup_iE_i)$,
where $\{E_i\}$ is the collection of all exceptional divisors. Equivalently
the isomorphism $f^{\sharp}$ is obtained  from the decomposition
$$\,\,g^*H^2(X;\ZZ)\op \op_i\ZZ c_1(E_i)=H^2(Z;\ZZ)
=h^*H^2(Y;\ZZ)\op \op_i\ZZ c_1(E_i)\,.$$
To prove the second statement we recall~[B, p.~772]  that Beauville's
quadratic form $B_X$ is the unique non-zero integral primitive
positive multiple  of the quadratic form $q_X$ defined by
$$\,\,q_X(\a)=
{n\over 2}\int_X\left(\tau_X\ov{\tau}_X\ri)^{n-1}\a^2+
(1-n)\int_X\tau_X^{n-1}\ov{\tau}_X^n\a\cdot
\int_X\tau_X^n\ov{\tau}_X^{n-1}\a\,.$$
Now express $q_X(\a)$ and $q_Y(\a)$ as integrals over $Z$ of the
appropriate forms pulled-back by $g$ and $h$ respectively, and use the
relation
$$\,\,g^*\a=h^*f^{\sharp}\a+\sum_in_iE_i\,,$$
 valid for some integers $n_i$. Then, as is easily checked, in order to
prove~(I.6.3) it suffices to show that
$$\left(g^*\tau_X g^*\ov{\tau}_X\ri)^{n-1}|_{E_i}=0=
\left(h^*\tau_Y h^*\ov{\tau}_Y\ri)^{n-1}|_{E_i}\eqno(\dag)$$
for any exceptional divisor $E_i$. Since $E_i$ is contracted by $g$ and by
$h$, and
since the fibers of $g|_{E_i}$ are distinct from the fibers of $h|_{E_i}$, at
every point of $E_i$ we have
$$\,\,\dim\hb{span of $\Ker(g|_{E_i})$ and $\Ker(g|_{E_i})$}\ge 2\,.$$
The above span is clearly contained in
$Ker(g^*\tau_X|_{E_i})=Ker(h^*\tau_Y|_{E_i})$, and hence this kernel is at
least two-dimensional.  This implies~($\dag$).
\qed
\msk

The proposition above allows us to identify the Beauville form of
$\cM_r^{2n}$ with that of $S^{[n]}$. The latter is described as follows~[B,
p.~777-778]. There is an inclusion $\s\cl
H^2(S;\ZZ)\hra H^2(S^{[n]};\ZZ)$ obtained composing the natural
"symmetrization map" $H^2(S;\ZZ)\hra H^2(S^{(n)};\ZZ)$, where $S^{(n)}$ is
the $n$-fold symmetric product of $S$, with the pull-back map $\e^*\cl
H^2(S^{(n)})\to H^2(S^{[n]})$, where $\e\cl S^{[n]}\to S^{(n)}$ is the
morphism mapping a subscheme to the $0$-cycle associated to it.  One has
$$\,\,H^2(S^{[n]};\ZZ)=\s\left(H^2(S;\ZZ)\ri)\op \ZZ T\,,\eqno({\rm I}.6.4)$$
where $T\in H^2(S^{[n]};\ZZ)$ is the (unique)  class such
that $2T$ is cohomologous to the divisor parametrizing
non-reduced subschemes.

\proclaim (I.6.5) Description of $B_{S^{[n]}}$.
The direct sum~(I.6.4) is orthogonal for $B_{S^{[n]}}$.  The restriction
of  $B_{S^{[n]}}$ to $\s\left(H^2(S;\ZZ)\ri)$ is equal to the intersection form
on $H^2(S;\ZZ)$, and $B_{S^{[n]}}(T)=-2(n-1)$.

Let's go back to the proof of Proposition~(I.0.4).  Since  the map $\t_r$
clearly preserves type, all we must  show is that the  homomorphism of
lattices is in fact an isomorphism. Since the Mukai form $\la,\ra$   has
discriminant one, the discriminant of its restriction to $v_r^{\bot}$ is
equal to $\la v_r,v_r\ra=2(n-1)$. In particular $\t_r$ is injective because
$n>1$. Since $\rk v_r^{\bot}=23=\rk H^2(S^{[n]})=\rk H^2(\cM_r^{2n})$ the
image $\t_r\left(H^*(S;\ZZ)\ri)$ has finite index, say $s$,  in
$H^2(\cM_r^{2n};\ZZ)$. By Proposition~(I.6.1) the map $\t_r$ is a
homomorphism of lattices, and hence  the discriminant of
$B_{\cM_r^{2n}}$ is equal to $s^2(2n-2)$.   By~(I.6.2) together
with~(I.6.5) the discriminant of $B_{\cM_r^{2n}}$ is equal to $2(n-1)$,
thus $s=1$, i.e.~$\t_r$ is an isomorphism.
\qed
\msk

\n
{\bf Proof of Proposition~(I.6.1).}
\hskip 2mm
Let $\cU_r\ss\cM_r$ be the open subset constructed in
Subsections~I.4-I.5. By~(I.4.2)-(I.4.31)-(I.5.17) the complement of $\cU_r$
in $\cM_r$ has codimension two, and hence $H^2(\cM_r)\cong H^2(\cU_r)$.
Thus $\t_r$ is determined by the map $\t_{\cF^r}\cl H^*(S)\to H^2(\cU_r)$
given by
$$\,\,\t_{\cF^r}(\a):=\rho_*\left[ch(\cF^r)^*(1+\o)\a\right]_3\,,$$
where $\rho\cl S\tm\cU_r\to\cU_r$ is projection, and
we have denoted by the same symbol  classes in $H^*(S)$ and their pull-back
to $H^*(S\tm\cU_r)$.  By~(I.4.31) and~(I.5.1) there is an isomorphism
$\vf\cl \cU_r\cong\cU_1$, and $\cU_1$ is an open subset
of $S^{[n]}$ whose complement has codimension two.  Thus by
Proposition~(I.6.2) $\vf$ induces an isomorphism of the lattices
$\left(H^2(\cM_r;\ZZ),B_{\cM_r}\ri)$ and
$\left(H^2(S^{[n]};\ZZ),B_{S^{[n]}}\ri)$. We will always tacitly
 identify  these lattices.

The proof of Proposition~(I.6.1) will be by induction on $r$.   By the
construction of $\cF^{r+1}$ we have two exact sequences:
$$\displaylines{\hfill\,\,0\to\cF^{r+1}\to\cE^{r+1}
\to\i_*\cQ_{r+1}\to 0\,,\hfill\llap{(I.6.6)}\cr
\hfill\,\,0\to \cO_S\tm\cU_r\to\cE^{r+1}\to
\cF^r\ot[-2C]\ot\rho^*\xi_{r+1}\to 0\,.\hfill\llap{({\rm I}.6.7)}\cr}$$
Here $\i$ is the inclusion of $S\tm(L\cup\Si_1)$ if $r=1$, and of
$S\tm\Si_1$ if $r\ge 2$, while $\cQ_{r+1}$ is a rank-one sheaf.
{}From the two  exact sequences above  one gets
$$\,\,\t_{\cF^{r+1}}(\a)=\t_{\cF^r}(\a\cdot e^{2C})
+\langle v_r,\a\cdot e^{2C}\rangle c_1(\xi_{r+1})
-\rho_*\left[ ch(\i_*\cQ_{r+1})^*(1+\o)\a\right]_3\,.\eqno({\rm I}.6.8)$$
Since $\xi_{r+1}=Ext^1_\rho(\cF^r\ot[-2C],\cO_{S\tm\cU_r})$,
using~(I.4.3)-(I.5.2)-(I.5.5) we get
$c_1(\xi_{r+1})=c_1(\rho_!(\cF^r\ot[-2C])$.
Applying Grothendieck-Riemann-Roch one gets
$$\,\,c_1(\xi_{r+1})=-\t_{\cF^r}(e^{2C})+
\rho_*\left[c_1(\cF^r)\o\right]_3\,.\eqno({\rm I}.6.9)$$
We need some notation. First we recall that if $\a\in H^*(S)$ then $\a^i$ is
the component of $\a$ belonging to $H^{2i}(S)$, and that $\o\in H^4(S;\ZZ)$
is the fundamental class. For $\a\in H^*(S)$ we set $\a^1_1:=\s(\a^1)\in
H^2(S^{[n]})$. If $\a^1$ is  integral and  $A\ss S$ is a real
surface representing it (i.e.~representing its Poincar\'e dual),   then
$\a^1_1$ is represented by the  subset of $S^{[n]}$ parametrizing schemes
intersecting $A$. By abuse of notation we will often denote by the same
symbol a divisor and its first Chern class.
\msk

\n
{\bf Rank one.}
\hskip 2mm
Since  $\cF^1=I_{\cZ}\ot[\Si+nC]$, we have
$$\,\,\t_{\cF^1}(\a)=
\rho_*\left[ch(I_{\cZ})^*e^{-\Si-nC}(1+\o)\a\right]_3\,.$$
First we compute $ch(I_{\cZ})$. Let $j\cl\cZ\hra S\tm\cU_1$ be the
inclusion, and let
$$\,\,\G:=\{(p,[Z])\in S\tm\cU_1|
\hb{ $p$ is a non-reduced point of $Z$}\}\,.$$
Then
$$\eqalign{ch(j_*\cO_{\cZ}) &
=[\cZ]-{1\over 2}[\G]+\hb{(higher order)},\cr
ch(I_{\cZ})=1-ch(j_*\cO_{\cZ})
&=1-[\cZ]+{1\over 2}[\G]+\hb{(higher order)},\cr}$$
where the first equality follows from Grothendieck-Riemann-Roch.
Substituting in the expression for $\t_{\cF^1}$ one gets
$$\t_{\cF^1}(\a)
=-\rho_*\left\{[\cZ]\cdot\left(\a^1-\a^0(\Si+nC)\right)\ri\}
-{1\over 2}\a^0\rho_*[\G]\,.$$
Since for $\b\in H^2(S)$ we have  $\rho_*\left\{[\cZ]\cdot\b\ri\}=\b_1$,
and since $\rho_*[\G]=2T$, we get
$$\,\,\t_{\cF^1}(\a)=-\a^1_1+\a^0(\Si_1+nC_1-T)\,.\eqno({\rm I}.6.10)$$
In particular $\t_{\cF^1}$ is integral, and hence also $\t_1$. Using
Proposition~(I.6.2) together with~(I.6.5) we conclude that
$$\,\,B_{\cM_1}\left(\t_{\cF^1}(\a),\t_{\cF^1}(\a)\right)=
\a^1\cdot\a^1-2\a^0\left(\a^1\cdot(\Si+nC)\ri)\,.$$
Since $\langle\a,v_1\rangle=0$ is equivalent to
$\a^2=\a^1\cdot(\Si+nC)$, the map $\t_1$ is a homomorphism of lattices.
This  finishes the proof of Proposition~(I.6.1) when the rank is one.
\msk

\n
{\bf Rank two.}
\hskip 2mm
Applying Equation~(I.6.9) one gets
$$\,\,c_1(\xi_2)=-(n-2)C_1-\Si_1+T\,.$$
Now let's compute
$\rho_*\left[ch(\i_*\cQ_2)^*(1+\o)\a\right]_3$. Let
$$\i^{\Si_1}\cl S\tm\Si_1\hra S\tm\cU_2,\qquad
\i^{L}\cl S\tm L\hra S\tm\cU_2$$
be the inclusions, and set $\cQ^{\Si_1}:=\cQ_2|_{S\tm\Si_1}$,
$\cQ^L:=\cQ_2|_{S\tm L}$. Then
$$\,\,\cQ^{\Si_1}\cong\pi^*[C]\ot\rho^*\l_1,\qquad
 \cQ^L\cong\pi^*[C]\ot I_{\O}\ot\rho^*\l_2\,,$$
where $\l_1$, $\l_2$ are line bundles on $\cU_2$, and
$$\,\,\O:=\left\{(P,x)|
\hb{ $x\in L$ and $P$ is the singular point
of $\cQ_2|_{S\tm\{x\}}$}\ri\}\,.$$
Applying Grothendieck-Riemann-Roch to $\i^{\Si_1}$, $\i^L$ one gets
$$\eqalignno{ch(\i^{\Si_1}_{*}\cQ^{\Si_1})&\equiv
\rho^*[\Si_1]+\pi^*[C]\cdot\rho^*[\Si_1]\pmod{H^3(\cU_2)},
&({\rm I}.6.11)\cr
ch(\i^L_{*}\cQ^L)&\equiv
\rho^*[L]+\pi^*[C]\cdot\rho^*[L]-[\i^L_{*}\O]\pmod{H^3(\cU_2)}.
&({\rm I}.6.12)}$$

We claim that the following linear equivalence among divisors on $S^{[n]}$
holds:
$$\,\,L\sim (n-1)C_1-T\,.\eqno({\rm I}.6.13)$$

\n
{\bf Proof of~(I.6.13).}
\hskip 2mm
Since $S^{[n]}$ is simply connected it suffices to prove that the divisors
are cohomologous, and for this it suffices to verify that they intersect
 any $2$-homology class in the same number of points. If $\b$ is a
two-cycle on $S$, and $W\ss S$ is a  set of $(n-1)$ points disjoint from
$\b$, let
$$\,\,\b_W:=\{[Z]\in S^{[n]}|\ Z=P\cup W \hb{, where } P\in\b\}\,.$$
Then $\b_W$ is a two-cycle on $S^{[n]}$, and   the annihilator of $c_1(T)$
is spanned by the $\b_W$. Clearly we have
$$\,\,\la c_1(L),\b_W\ra=(n-1)\la c_1(C),\b\ra=
(n-1)\la c_1(C_1)-c_1(T),\b_W\ra\,.$$
Now choose a point $P\in\Si$ and a subset $W\ss S$ of $(n-2)$ points
disjoint from $\Si$. Let $\wt{\Si}\ss S^{[n]}$ be the curve given by
$$\,\,\wt{\Si}:=\{[Z_Q]\in S^{[n]}|
\ Z_Q=Q\cup P\cup W \hb { where } Q\in\Si\}\,,$$
where the scheme structure of $Z_P$ at $P$ consists of the double
point contained in $\Si$. Since
$$\,\,\la c_1(L),\wt{\Si}\ra=n-2\,,\qquad
\la c_1(C_1),\wt{\Si}\ra=1\,,\qquad
\la c_1(T),\wt{\Si}\ra=1\,,$$
the intersections of the two sides of Equation~(I.6.13) with $\wt{\Si}$
are equal. Since $H_2(S^{[n]})$ is spanned by $\wt{\Si}$ and the classes
 $\b_W$, we conclude that~(I.6.13) holds.
\qed
\msk

Using~(I.6.13) together with Equations~(I.6.11)-(I.6.12) one  obtains that
$$\displaylines{\quad\rho_*\left[ch(\i_*\cQ_2)^*(1+\o)\a\right]_3=
(n-1)\left(\int_S \a^1\cdot C-\a^2\ri)C_1\hfill\cr
\hfill{}-\left(\int_S \a^0\o-\a^1\cdot C+\a^2\ri)\Si_1
-\left(\int_S \a^1\cdot C-\a^2\ri)T\,.\quad\cr}$$
At this point we have all the elements needed to apply Formula~(I.6.8).
The result is
$$\displaylines{\quad\t_{\cF^2}(\a) = -\a^1_1
-\left(\int_S(n-2)\a^0\o+(n-2)\a^1\cdot\Si
+(n^2-3n+3)\a^1\cdot C-(2n-3)\a^2\right)C_1 \hfill\cr
\hfill{}-\left(\int_S\a^1\cdot\Si+(n-1)\a^1\cdot C-2\a^2\right)\Si_1
+\left(\int_S\a^0\o+\a^1\cdot\Si
+(n-1)\a^1\cdot C-2\a^2\right)T\,.\quad\cr}$$
In particular $\t_{\cF^2}$ is integral.
Since $\langle\a,v_2\rangle=0$ is equivalent to
$\a^0=\left(2\a^2-\a^1\cdot\Si-(n-2)\a^1\cdot C\right)$, one gets that
$$\,\,\t_2(\a)=-\a_1^1-\left(\int_S (n-1)\a^1\cdot C-\a^2\ri)C_1
-\left(\int_S\a^1\cdot\Si+(n-1)\a^1\cdot C-2\a^2\ri)\Si_1
+\left(\int_S\a^1\cdot C\ri)T\,.$$
A tedious but straightforward computation shows that $\t_2$ is a
homomorphism of lattices.
\msk

\n
{\bf Rank at least three.}
\hskip 2mm
Let $\O\ss H^*(S;\ZZ)$ be the $\ZZ$-span of $\{1,C,\Si,\o\}$.
Since the restriction of the Mukai form to $\O$ is unimodular,
 $H^*(S;\ZZ)$ is the orthogonal direct sum of $\O$ and of $\O^{\bot}$.
Notice that $\O^{\bot}\ss v_r^{\bot}$, and thus
$$\,\,v_r^{\bot}=
\left(v_r^{\bot}\cap\O\ri)\op_{\bot}\O^{\bot}\,.\eqno({\rm I}.6.14)$$
Let $r\ge 1$. By Formula~(I.6.8)  the restrictions of  $\t_{\cF^{r+1}}$
and $\t_{\cF^r}$ to $\O^{\bot}$ are  equal, and hence by induction we have
$$\t_r|_{\O^{\bot}}=\t_1|_{\O^{\bot}}\eqno({\rm I}.6.15)$$
for all $r\ge 1$. In order to examine the restriction of $\t_r$ to
$v_r^{\bot}\cap\O$, it is convenient to introduce the class
$\b_r:=v_r-2(n-1)C$; notice that $\b_r\in v_r^{\bot}\cap\O$. A computation
gives
$$\,\,\la\b_r,\b_r\ra=-2(n-1)\,.\eqno({\rm I}.6.16)$$
Let $\L_r:=\{v_r,\b_r\}^{\bot}\cap\O$. As is easily checked we have an
orthogonal direct sum decomposition
$$\,\,v_r^{\bot}\cap\O=\ZZ\b_r\op_{\bot}\L_r\,.\eqno({\rm I}.6.17)$$
A computation gives
$$\,\,\L_r:=\{x+yC+z\o|\  (r-1)x+y-rz=0\}\,.\eqno({\rm I}.6.18)$$

\proclaim (I.6.19) Lemma.
Keep notation as above. If $r\ge 3$ then
$$\,\,\t_{\cF^r}(x+yC+z\o)=x\left(C_1-(r-2)\Si_1\ri)-y\Si_1+z(r\Si_1-C_1),
\qquad \t_{\cF^r}(\b_r)=T\,.$$

\n
We will prove the lemma at the end of the section. First we finish the
proof  of Proposition~(I.6.1)  assuming the lemma. By
Equality~(I.6.17) and Lemma~(I.6.19) the restriction of $\t_r$ to
$v_r^{\bot}\cap\O$ is integral. By~(I.6.15)  the restriction of $\t_r$
to $\O^{\bot}$ is also integral. By~(I.6.14) we conclude that $\t_r$ is
integral.   Now let's prove that $\t_r$ is a homomorphism of lattices.
Since $\t_1$ is a lattice homomorphism, and since $\t_1(\O^{\bot})$ is the
orthogonal complement of $\{C_1,\Si_1,T\}$, it follows  that
    $\t_r|_{\O^{\bot}}$ is a lattice homomorphism (use~(I.6.15)) and that
$\t_r(\O^{\bot})$ is perpendicular to $\{C_1,\Si_1,T\}$. By
Lemma~(I.6.19) the image $\t_r(v_r^{\bot}\cap\O)$ is contained in the
$\ZZ$-span of $\{C_1,\Si_1,T\}$, hence it is perpendicular to
$\t_r(\O^{\bot})$. It remains to verify  that the restriction of
$\t_r$ to $v_r^{\bot}\cap\O$ is a lattice homomorphism.  If we restrict
$\t_r$ to $\ZZ\b_r$ we get a lattice homomorphism by Formula~(I.6.16).
Lemma~(I.6.19) shows that $\t_r(\b_r)$ is perpendicular to
$\t_r(\L_r)$. Hence from~(I.6.17) we see that we are
reduced to proving that the restriction of $\t_r$ to $\L_r$ is a lattice
homomorphism. This consists of a straightfroward computation
(use~(I.6.18)).
\msk

\n
{\bf Proof of Lemma~(I.6.19).}
\hskip 2mm
It follows from Exact sequences~(I.6.6)-(I.6.7) that for $r\ge 2$ one has
$$\,\,\rho_*\left[c_1(\cF^{r+1})\o\right]_3  =
\rho_*\left[c_1(\cF^{r})\o\right]_3+rc_1(\xi_{r+1})-\Si_1\,.$$
Equation~(I.6.8) gives   the following formulae,  for $r\ge 2$:
$$\eqalign{
\t_{\cF^{r+1}}(1) & = \t_{\cF^r}(e^{2C})+(r+1)c_1(\xi_{r+1})+\Si_1\,,\cr
\t_{\cF^{r+1}}(C) & = \t_{\cF^r}(C)+c_1(\xi_{r+1})\,,\cr
\t_{\cF^{r+1}}(\Si) & =
\t_{\cF^r}(\Si+2\o)-(r^2+r+2-n)c_1(\xi_{r+1})-\Si_1\,,\cr
\t_{\cF^{r+1}}(\o) & =
\t_{\cF^r}(\o)-rc_1(\xi_{r+1})+\Si_1\,.\cr}\eqno({\rm I}.6.20)$$
Using~(I.6.6)-(I.6.7) and~(I.6.13)     one gets
$$\,\,\rho_*\left[c_1(\cF^2)\o\right]_3 = -(2n-3)C_1-2\Si_1+2T\,.$$
Formula~(I.6.9)  together with the equality above and the computations
for rank two give
$$\,\,c_1(\xi_3)=(n-1)C_1-T\,.$$
Plugging  this value in Formulae~(I.6.20), and
using the computations for rank two,  one obtains Lemma~(I.6.19) for $r=3$.
In order to prove Lemma~(I.6.19) for $r\ge 4$ we add to the lemma the
following statement.

\proclaim (I.6.21).
Keep notation as above. If $r\ge 4$ then
$$\eqalign{\rho_*\left[c_1(\cF^{r-1})\o\right]_3 & =-r\Si_1+C_1\,,\cr
c_1(\xi_r) & =0\,.\cr}$$

The proof of Lemma~(I.6.19) for $r\ge 4$ is by induction.
First one verifies that~(I.6.21) and Lemma~(I.6.19) hold for $r=4$ by a
computation, and then the inductive step is an easy consequence of
Formulae~(I.6.20).
\bsk

\n
{\bf  II. Proof of the main theorem.}
\msk
\n
We will   need the following technical result; its proof is deferred to the
Section~IV.

\proclaim (II.1) Proposition.
Let $S$ be a $K3$ surface, and $v\in H^*(S;\ZZ)$ be a Mukai vector. Suppose
that $v^0$ and the order of divisibility of $v^1$ are coprime.
\msk
\item{1.} If an ample divisor on $S$  is  $|v|$-generic  then it is also
$v$-stabilizing.
\item{2.} Let $H$ be a $v$-stabilizing ample divisor on $S$ (for example
this is the case if $H$ is $|v|$-generic, by Item~(1)). Let
$L$ be a $|v|$-generic ample divisor on $S$, and let $\cC$ be the unique
open $|v|$-chamber containing $L$. If $H\in\ov{\cC}$ then
$H$-slope-stability is the same as $L$-slope-stability, and hence
$\cM_v(H)=\cM_v(L)$.

By the above proposition we can assume, in proving the Main Theorem, that
$H$ is $v$-generic.   We will prove Theorem~(2)  by first deforming $S$ to
an elliptic surface, and then invoking Theorem~(I.0.4). Let $\cK_{2d}$ be the
moduli space of polarized K3 surfaces of degree $2d$, i.e.~couples
$(S,H)$, where $S$ is a $K3$ and $H$ is a primitive ample divisor
on $S$ with $H^2=2d$. Let $[S,H]\in K_{2d}$, and assume $v\in H^*(S;\ZZ)$
is a Mukai vector such that
$$\,\,v^1=\pm c_1(H)\,.\eqno({\rm II}.2)$$
Given any  $[S_x,H_x]\in\cK_{2d}$   the Mukai vector
$v$ makes sense in $H^*(S_x;\ZZ)$, because of~(II.2), and hence we can
consider the moduli space $\cM_v(S_x,H_x)$. Let
$$\,\,\cK_{2d}(v):=
\{[S_x,H_x]\in\cK_{2d}| \hb{ $H_x$ is $v$-generic}\}\,.$$
This is an open subset of $\cK_{2d}$; in fact its complement is the union of
a finite set of components of the Noether-Lefschetz locus.

\proclaim (II.3) Proposition.
Keep notation as above, and assume that~(II.2) is satisfied. Suppose
there  exists $[T,L]\in\cK_{2d}(v)$ such that Theorem~(2) holds for the
moduli space $\cM_v(T,L)$. Given any other $[S,H]\in\cK_{2d}(v)$
Theorem~(2) holds for the moduli space $\cM_v(S,H)$.

\pf
Let $\cH_{2d}(v)$ be the parameter space (open subset of a Hilbert scheme)
for $K3$ surfaces $[S,H]\in\cK_{2d}$ embedded in  projective space by
a high multiple of $H$. Let $\cS_{2d}\to\cH_{2d}(v)$ be the
tautological family of (embedded) surfaces. By~[Ma] there exists a relative
moduli space $\pi\cl\cM_v(2d)\to\cH_{2d}(v)$: the fiber of $\pi$ over a
surface $S_x$ embedded by a multiple of $H_x$ is isomorphic to
$\cM_v(S_x,H_x)$. By~[M1, Th.~(1.17)] $\pi$ is a submersion at every
point (here we use the fact that on $\cK_{2d}(v)$ semistability implies
stability). Since $\cH_{2d}(v)$ is irreducible (by irreducibilty of
$\cK_{2d}(v)$) and $\pi$ is
proper~[Ma], the map   $\pi$ is surjective. Hence, by irreducibilty of
$\cH_{2d}(v)$, the moduli space $\cM_{v}(S,H)$
is a deformation  of $\cM_v(T,L)$, and therefore Item~(1) of Theorem~(2)
holds for $\cM_v(S,H)$. To prove Item~(2) notice that one can
construct a relative quasi-tautological family of sheaves on
$\cS_{2d}\tm_{\cH_{2d}(v)}\cM_v(2d)$ (see the proof of [M2, Th.~(A.5)]).
Hence the map $\t_{v}\cl v^{\bot}\cap H^*(S_x)\to H^2(\cM_v(S_x,H_x))$ is a
locally constant function of $x\in\cH_{2d}(v)$. Since $\cH_{2d}(v)$ is
irreducible we conclude that Item~(2) of Theorem~(2) holds for
$\cM_v(S,H)$.
\qed
\msk

First we prove Theorem~(2) under particularly favourable hypotheses.

\proclaim (II.4) Proposition.
Let hypotheses be as in the statement of Theorem~(2). In addition suppose
that $v^1=c_1(H)$, where $H$ is a (primitive) ample divisor on $S$, and that
$$\,\,H^2\ge 2|v|\,,\qquad H^2\ge 4\,.\eqno({\rm II}.5)$$
Then Theorem~(2) holds for the moduli space $\cM_v(S,H)$.

\pf
Let $T$ be an elliptic $K3$ with $Pic(S)=\ZZ[\Si]\op\ZZ[C]$, where $\Si$ is
a section of the elliptic fibration, and $C$ is an elliptic fiber. Let
$d:=H^2/2$, and set $L:=(\Si+(d+1)C)$. By~(II.5) we have $d\ge 2$,
hence $L$ is an ample divisor on $T$.  Since $L$ is primitive, and since
$L^2=2d=H^2$, both $[S,H]$ and $[T,L]$ belong to the same moduli space
$\cK_{2d}$. By~(II.5) together with Lemma~(I.0.3) the polarization $L$ is
$|v|$-suitable, in particular $|v|$-generic. The polarization $H$ is
$|v|$-generic by hypotheses, hence $[T,L],[S,H]\in\cK_{2d}(v)$.
Since $L$ is $|v|$-suitable, and since $v^1$ ($=c_1(L)$) is a numerical
section,  Theorem~(I.0.4) tells us that the Main Theorem holds for
$\cM_v(T,L)$. By Proposition~(II.3) we conclude that the Main Theorem
holds also for $\cM_v(S,H)$.
\qed
\msk

Next we concentrate on the case when $\rho(S)$ (the rank of $Pic(S)$) is at
least two, with no further hypotheses.

\proclaim (II.6) Lemma.
Let $S$ be a a projective $K3$ surface with $\rho(S)\ge 2$.
Let $v\in H^*(S;\ZZ)$ be a Mukai
vector with $v^1$ primitive, and $\cC$ be  an open $v$-chamber.     There
exists a Mukai vector $w$ equivalent to $v$ such that $w^1=c_1(H)$, where
$H$ is  a primitive ample divisor belonging to $\cC$.
We can choose  $w$ so that $H^2$ is arbitrarily large.

\pf
Let's show that there exists an ample divisor $H\in\cC$ such
that
\msk
\item{a.} $H$ is primitive, and
\item{b.} $c_1(H)\equiv v^1\pmod{v^0H^{1,1}_{\ZZ}}$.
\msk
\n
There exists a $\ZZ$-basis of $H^{1,1}_{\ZZ}(S)$ consisting of elements
of $\cC$; we identify $H^{1,1}_{\ZZ}(S)$ with $\ZZ^{\rho}$ via such a
basis. There exist $\ul{\rm positive}$ integers $m$, $a_i$ (for
$i=1,\ldots,\rho$), with $gcd\{a_1,\ldots,a_{\rho}\}=1$,  such that
$$\,\,v^1\equiv (ma_1,\ldots,ma_{\rho}) \pmod{v^0H^{1,1}_{\ZZ}}\,.$$
Since $v^1$ is primitive, we have
$gcd\{m,v^0\}=1$. Let $H$ be the divisor (class) such that
$$\,\,c_1(H)=
(ma_1+v^0n_1,\ldots,ma_{\rho-1}+v^0n_{\rho-1},ma_{\rho})\,,$$
where,  for all $i=1,\ldots,\rho-1$,
\msk
\item{$\bu$} $n_i>0$,
\item{$\bu$} $gcd\{m,n_i\}=1$, and
\item{$\bu$} if $p$ is a prime dividing $a_{\rho}$ but not $m$, then $p$
divides each $n_i$.
\msk
\n
Since all the coordinates of $c_1(H)$ are strictly positive, $H$
belongs to $\cC$, in particular it is ample.  By construction
$H$ satisfies Items~(a) and~(b). If the $n_i$ are arbitrarily large then
$H^2$ is arbitrarily large.
 Now let $\xi$ be the line bundle on $S$ such that
$$c_1(H)=v^1+v^0c_1(\xi)\,,\quad
\hb{ i.e.~$c_1(\xi)=(n_1,\ldots,n_{\rho-1},0)$.}$$
If $w:=ch(\xi)v$ then $w^1=c_1(H)$, hence $w$ satisfies
the conclusions of the lemma.
\qed

\proclaim (II.7) Proposition.
Let hypotheses be as in the statement of Theorem~(2). In addition
suppose  that $\rho(S)\ge 2$. Then Theorem~(2) holds for $\cM_v(H)$.

\pf
Let $\cC$ be an open $|v|$-chamber such that
$H\in\ov{\cC}$. By Lemma~(II.6) there exists a Mukai vector $w$ equivalent
to $v$ such that $w^1=c_1(L)$, where $L$ is a primitive ample divisor
such that $L\in\cC$ and
$$\,\,L^2\ge 2|v|\,,\qquad L^2\ge 4\,.$$
Since $w$ is equivalent to $v$ we have $\cM_w(L)=\cM_v(L)$, and by
Proposition~(II.1) $\cM_v(L)=\cM_v(H)$.  Hence it it is equivalent  to
prove that Theorem~(2) holds for $\cM_w(L)$. By the above inequality we
can apply Proposition~(II.4), hence the Main Theorem does indeed hold for
$\cM_w(L)$.
\qed
\msk

We are left with proving that the Main Theorem holds when $\rho(S)=1$.
Assume $S$ is such a surface.  Let $H$ be the ample generator of $Pic(S)$,
and set  $d:=H^2/2$. Since $v^1$ is primitive we have $v^1=\pm c_1(H)$,
hence we can apply  Proposition~(II.3). By surjectivity of the period map for
polarized $K3$ surfaces there exists  $[T,L]\in\cK_{2d}(v)$ with
$\rho(T)\ge 2$. By Proposition~(II.7) the Main Theorem holds
$\cM_v(T,L)$, and thus Proposition~(II.3) shows that   Theorem~(2) holds
for $\cM_v(S,H)$. This finishes the proof of the Main Theorem.
\bsk

\n
{\bf  III. Higher-rank Donaldson polynomials.}
\msk
\n
Let $(S,H)$ be a polarized $K3$ surface, and $v\in H^*(S;\ZZ)$ be a Mukai
vector; we assume that $H$ is $v$-stabilizing. If the rank (i.~e.~$v^0$) is
two then one can associate to $\cM_v$ a certain Donaldson polynomial
function on $H_2(S)$. The (algebro-geometric) definition extends to any
rank as follows. Let $\cF$ be a quasi-tautological family of sheaves on
$S\tm\cM_v$, and let $\rho\cl S\tm\cM_v\to\cM_v$ be projection.
We define $\mu_{\cF}\cl H_2(S)\to H^2(\cM_v)$ by setting
$$\,\,\mu_{\cF}(\a):=
\rho_*\left[\left(-{1\over \s(\cF)}ch_2(\cF)+
{1\over 2\s(\cF)^2\cdot\rk(\cF)}ch_1^2(\cF)\ri)\cdot\a\ri]\,.$$
If $\cG$ is another quasi-tautological family of sheaves, then  there exist
vector bundles $\xi$, $\eta$ on $\cM_v$ such that
$\cF\ot\rho^*\xi\cong\cG\ot\rho^*\eta$~[M2, Th.~(A.5)]; it follows
that $\mu_{\cF}=\mu_{\cG}$, and hence we can set
$$\mu_v:=
\mu_{\cF}\qquad \hb{ where $\cF$ is any quasi-tautological family.}$$
Now let $q_v\cl H_2(S;\QQ)\to \QQ$ be the polynomial given by
$$\,\,q_v(\a):=\int_{\cM_v}\mu_v(\a)^{d(v)}\,,$$
where $d(v):=2+\la v,v\ra$ is the (complex) dimension of $\cM_v$. If
$v^0=2$ then $q_v$ equals the corresponding Donaldson polynomial by a
theorem of Morgan~[Mo]. Donaldson's polynomials of $K3$ surfaces have been
computed long ago; we will show that the formula for rank two holds also in
higher rank.

\proclaim (III.1) Proposition.
Keeping notation as above, assume that $v^1$ is primitive. Set $d(v)=2n$,
and let $Q\cl H_2(S;\QQ)\to\QQ$ be the intersection form. Then
$$\,\,q_v={(2n)!\over n!2^n}Q^n\,.$$

The above result is a straightforward corollary of our Main
Theorem~(2) if $\dim\cM_v>2$, and of results of Mukai if $\dim\cM_v\le
2$, together with a theorem of Fujiki~[Fu, Th.~(4.7)]. We will provide a
simple proof of Fujiki's result avoiding hyperk\"ahler structures.

\proclaim Theorem (Fujiki~[Fu, Th.~(4.7)]).
Let $X$ be an irreducible K\"ahler symplectic compact manifold of
(complex) dimension $2n$. There exists $\l_X\in\QQ$ (call it the {\it Fujiki
constant of X}) such that
$$\int_X\b^{2n}=\l_X\cdot B_X(\b)^n \eqno({\rm III}.2)$$
for all $\b\in H^2(X)$, where $B_X\cl H^2(X)\to \CC$ is the quadratic form
associated to $X$ by Beauville~[B, Th.~(5)]. Fujiki's constant is
independent of the birational (symplectic) model of $X$.

\pf
Let $p_X\cl H^2(X)\to\CC$ be the degree-$2n$ polynomial defined by the
left-hand side of~(III.2). Consider the versal deformation space $f\cl\cX\to
D$ of $X$, and for $s\in D$ let $X_s:=f^{-1}(s)$. Let $u_s\cl
H^2(X_s)\brel\cong\over\to H^2(X)$ be the isomorphism given by the
Gauss-Manin connection. Then
$$\,\,p_X(u_s\b)=p_{X_s}(\b)\,.\eqno(*)$$
If $\vf_s\in H^{2,0}(X_s)$ then $p_{X_s}(\vf_s)=0$ by type consideration,
and  thus by~($*$) we have
$$\,\,p_X(u_s\vf_s)=0\,.$$
By Beauville~[B, Th.~(5)] the collection
$$\{[u_s\vf_s]|\ s\in D,\ 0\not=\vf_s\in H^{2,0}(X_s)\}$$
fills out an open (in the analytic topology) subset of the smooth quadric
where $B_X$ vanishes. Hence we deduce that $p_X=B_X\cdot p_X'$
for some degree-$2(n-1)$ polynomial $p_X'$. Now consider
$$\,\,p_X(\ub{u_s\vf_s,\ldots,u_s\vf_s}_{2n-1},u_s\ov{\vf_s})=
n^{-1}B_X(u_s\vf_s,u_s\ov{\vf_s})\cdot p_X'(u_s\vf_s)\,.$$
Since the left-hand side vanishes by type consideration (unless $n=1$, in
which case we are done), and since $B_X(u_s\vf_s,u_s\ov{\vf_s})>0$
by~[B, Th.~5], we conclude that $p_X'(u_s\vf_s)=0$ for all $\vf_s$, and
hence $B_X|p_X'$. Thus $p_X=B_X^2p_X''$ for some
degree-$2(n-2)$ polynomial $p_X''$. Evaluating
$$\,\,p_X(\ub{u_s\vf_s,\ldots,u_s\vf_s}_{2n-i},
\ub{u_s\ov{\vf_s},\ldots,u_s\ov{\vf_s}}_i)
\qquad i\le (n-1)\,,$$
 and arguing similarly
we conclude that $p_X=\l_XB_X^n$ for some constant $\l_X$. That
$\l_X\in\QQ$ follows from integrality of $B_X$~[B] and of $p_X$. Finally
let's show that if $f\cl X\cdots>Y$ is birational then $\l_X=\l_Y$. Since
$f,f^{-1}$ are isomorphisms in codimension one they induce an isomorphism
$f^{\sharp}\cl H^2(X)\brel\cong\over\to H^2(Y)$, and one has
$B_X(\b)=B_Y(f^{\sharp}\b)$ for all $\b\in H^2(X)$. We claim that if
$\vf_X\in H^{2,0}(X)$ then
$$\,\,\l_X B_X(\vf_X+\ov{\vf_X})^n=\int_X(\vf_X+\ov{\vf_X})^{2n}
=\int_Y(f^{\sharp}\vf_X+f^{\sharp}\ov{\vf_X})^{2n}
=\l_Y B_Y(f^{\sharp}\vf_X+f^{\sharp}\ov{\vf_X})^n\,.$$
In fact let $X\supset U\cong V\subset Y$ be the
open subsets identified by $f$.  By Hartog's
Theorem $\vf_X|U=\vf_X|V$ extends to a holomorphic two-form on all of
$Y$, which necessarily represents $f^{\sharp}\vf_X$. Since the integrals
appearing above  can be computed  by integrating over $U$ and
$V$ respectively, we conclude that they  are equal. Since
$B_X(\vf_X+\ov{\vf_X})\not=0$ when $\vf_X\not=0$ we conclude that
$\l_X=\l_Y$.
\qed

\proclaim (III.3) Claim.
Let notation and hypotheses be as in Proposition~(III.1). The Fujiki constant
of $\cM_v$ is equal to $(2n!)!/n!2^n$.

\pf
If $\dim\cM_v=2$ then by a theorem of Mukai~[M2] the moduli space is a
$K3$ surface, and hence the claim holds. If $\dim\cM_v>2$ then by
Theorem~(2) the moduli space is deformation equivalent to a symplectic
birational model of $T^{[n]}$, call it $X$, where $T$ is another projective
$K3$, and $2n=d(v)$. In fact the proof of Theorem~(2) shows that $\cM_v$
can be deformed to $X$ through symplectic projective varieties. The Fujiki
constant is clearly   invariant under deformations through
symplectic K\"ahler  varieties, hence $\l_{\cM_v}=\l_X$. By Fujiki's
Theorem $\l_X=\l_{T^{[n]}}$. A simple calculation shows that
$$\,\,\l_{T^{[n]}}={(2n)!\over n!2^n}\,.$$
 This finishes the proof of the claim.
\qed
\msk

\n
{\bf Proof of Proposition~(III.1).}
\hskip 2mm
Assume first that $\dim\cM_v>2$. Let $\i_v\cl H^2(S)\to H^*(S)$ be the map
defined by
$$\,\,\i_v(\a):=\a+\left({1\over v^0}\int_Sv^1\wedge\a\ri)\o\,,$$
where $\o\in H^4(S;\ZZ)$ is the fundamental class. The image of $\i_v$ is
contained in $v^{\bot}$, hence it makes sense to compose $\i_v$ and $\t_v$.
A straightforward computation gives
$$\,\,\mu_v=-\t_v\circ\i_v\,.$$
If $Q$ is the intersection form, then $\i_v^*\la\a,\a\ra=Q(\a)$. Since
$\t_v$ is an isometry between $\left(v^{\bot},\la,\ra\ri)$ and
$\left(H^2(\cM_v),B_{\cM_v}\ri)$ (by the Main Theorem) we conclude that
$$\,\,B_{\cM_v}(\mu_v(\a))=Q(\a)\,.$$
The proposition follows at once from this equality together with Fujiki's
Theorem and Claim~(III.3).
Now assume that $\dim\cM_v=2$: the same argument as above applies
thanks to Mukai's description of $\left( H^2(\cM_v),B_{\cM_v}\ri)$ as
$v^{\bot}/\CC v$~[M2]. Finally when $\dim\cM_v=0$, the proposition is
equivalent to the statement that $\cM_v$ consists of a single reduced
point. We proved in Section~I that moduli spaces on an elliptic
surface, with $v^1$ a numerical section, are non-empty.   This implies that
the moduli space is a reduced point by~[M2, Cor.~(3.6)].
\bsk

\n
{\bf  IV. Appendix.}
\msk
\n
In this section we will prove Propositions~(II.1)-(I.1.6). We let $S$ be a
projective smooth irreducible surface. For a torsion-free sheaf $F$ on $S$
we set
$$\,\,\D_F:=c_2(F)-{r_F-1\over 2r_F}c_1(F)^2\,,$$
 where $r_F$ is the rank of $F$. We will be exclusively
interested in slope-(semi)stability, and hence  we will sistematically
omit the prefix "slope".

\proclaim  Lemma.
Let $H$ be a polarization of $S$. Let $F$ be a strictly
$H$-semistable torsion-free sheaf on $S$, and let
$$0\to A\to F\brel\pi\over\to B\to 0$$
be destabilizing, i.~e.~$\left[r_Fc_1(A)-r_Ac_1(F)\ri]\cdot H=0$.
Then
$$\,\,-{r_F^3\over 2}\D_F\le\left[r_Fc_1(A)-r_Ac_1(F)\ri]^2
\le 0\,.\eqno({\rm IV}.1)$$
Furthermore the right inequality is an equality only if
$\left[r_Fc_1(A)-r_Ac_1(F)\ri]=0$.

\pf
Replacing $A$ by $\pi^{-1}(Tor(B))$, where $Tor(B)$ is the torsion subsheaf
of $B$, we can assume $B$ is torsion-free.  First suppose $c_1(F)=0$. Since
$A$ is slope destabilizing $c_1(A)\cdot H=0$, and hence by Hodge index
$c_1(A)^2\le 0$, with equality only if $c_1(A)=0$. We are left with proving
the other inequality. Since $c_1(B)=-c_1(A)$ we have
$$\,\,c_2(F)=c_2(A)+c_2(B)-c_1(A)^2\,.$$
Both $A$ and $B$ are semistable, and thus by Bogomolov' inequality
$$\,\,c_2(A)\ge {r_A-1\over 2r_A}c_1(A)^2\qquad
c_2(B)\ge {r_B-1\over 2r_B}c_1(A)^2\,.$$
Substituting into the previous equality one gets the left inequality
in~(IV.1). If $c_1(F)\not=0$, we formally set $F_0:=F\ot\left(\det
F\ri)^{-1/r_F}$. Then $c_1(F_0)=0$, and applying the previous argument to
$F_0$ one gets~(IV.1).
\qed

\msk
\n
{\bf Proof of Item~(1) of Proposition~(II.1).}
\hskip 2mm
If $S$ is a $K3$, then
$$\,\,{r_F^3\over 2}\D_F=|v(F)|\,,$$
and hence, using the notation of the previous lemma, we have
$$\,\,-|v(F)|\le\left[r_Fc_1(A)-r_Ac_1(F)\ri]^2\le 0\,.$$
By the hypothesis on the divisibility of $c_1(F)$ we see that
$\left[r_Fc_1(A)-r_Ac_1(F)\ri]\not=0$, and thus
$\left[r_Fc_1(A)-r_Ac_1(F)\ri]^{\bot}\cap A(S)$ is a $v(F)$-wall. Since
$\left[r_Fc_1(A)-r_Ac_1(F)\ri]\cdot H=0$, $H$ belongs to this wall.

\proclaim (IV.2) Remark.
{\rm In fact we have proved something slightly stronger, namely that if $F$
is an $H$-semistable torsion-free sheaf on $S$ with $v(F)=v$, then $F$ is
$H$-stable.}

\msk
\n
{\bf Proof of Item~(2) Proposition~(II.1).}
\hskip 2mm
Item~(2) follows at once from the following.

\proclaim (IV.3) Lemma.
Let $S$ be a $K3$ surface, and let  $H_0$, $H_1$ be ample
divisors on $S$. Assume $F$ is a torsion-free sheaf on $S$ such  that the
the order of divisibility of $c_1(F)$ is coprime to $\rk(F)$.
If $F$ is $H_0$-slope-stable and $H_1$-slope-unstable, then there exists
a $|v(F)|$-wall $\xi^{\bot}$ such that $\xi\cdot H_0$, $\xi\cdot H_1$ are
of opposite signes (and non-zero ).

\pf
Since $|v(F^{**})|\le |v(F)|$ we can replace $F$ by $F^{**}$, and   assume
that $F$ is locally-free.  One can extend in the obvious way the notion of
(semi)stability to arbitrary elements  of $A(S)$.  Letting
$\L:=\RR_{+} c_1(H_0)\op \RR_{+} c_1(H_1)$, we set
$$\eqalign{\L^s:= & \{H\in \L|\hb{ $F$ is $H$-stable}\} \cr
\L^u:= & \{H\in \L|\hb{ $F$ is $H$-unstable}\}.} $$
Let $H\in\L^u$, and let $E$ be an $H$-desemistabilizing subsheaf of $F$;
then, by continuity,  $E$ is $H'$-desemistabilizing for all $H'$ near $H$, and
hence $\L^u$ is an open subset of $\L$. Let's show that also $\L^s$ is open
in $\L$. Set $H_t:=(1-t)H_0+tH_1$ for $t\in [0,1]$. Since $\L^s$ is a
(positive) cone it suffices to show that the set of $t$ such that
$H_t\in\L^s$ is open. Assume $H_{\ov{t}}\in\L^s$. Let $I\ss [0,1]$ be the
subset of points $t$ such that
$$\left[r_Fc_1(E)-r_Ec_1(F)\ri]\cdot H_t=0$$
for some  $H_1$-destabilizing subsheaf $E\ss F$ with torsion-free
quotient.  The set of $H_1$-destabilizing subsheaves of $F$
with torsion-free quotient is a bounded set (Grothendieck), and hence $I$ is
a finite set. Define $t_0$ by setting $t_0:=\min I$. Since $\L^s$ is a
convex cone,  we have $t>\ov{t}$ for all $t\in I$, and hence $t_0>\ov{t}$.
We claim that $F$ is $H_t$-stable for all $t<t_0$. In fact let $E\ss F$ be a
subsheaf with torsion-free quotient. If $E$ is not $H_1$-destabilizing then,
since $F$ is $H_0$-stable, $E$ is not $H_t$-destabilizing for all $t<1$, in
particular for $t<t_0$.  If $E$ is $H_1$-desemistabilizing then, since
$$\vf_E(t):=\left[r_Fc_1(E)-r_Ec_1(F)\ri]\cdot H_t$$
is a linear function of $t$, and since $\vf_E(1)>0$, $\vf_E(0)<0$, it will be
negative for all $t<t_E$, where $t_E\in[0,1]$ is the unique  solution of
$\vf_E(t)=0$. Clearly $t_E$ is a rational number and hence
$t_E\in I$. Thus $t_0\le t_E$, and we conclude that $E$ is not
desemistabilizing for all $t<t_0$. This finishes the proof that $\L^s$ is
open. Now we can prove Lemma~(IV.3). The nonempty
subsets $\L^s,\L^u\ss\L$ are open and disjoint; since $\L$ is connected
there must exist $H\in\L$ not belonging to either of these subsets,
i.~e.~$F$ is strictly $H$-semistable. Replacing $H$ by  an appropriate
multiple of itself, we can assume $H$ is an integral  class.  By
Item~(1) of Proposition~(II.1) we conclude that $H$ belongs to a
$|v(F)|$-wall; this wall separates $H_0$ from $H_1$.
\msk

\n
{\bf Proof of Proposition~(I.1.6).}
\hskip 2mm
By definition a $|v|$-suitable polarization is also $|v|$-generic, and
hence Item~(1) follows from  Proposition~(II.1). In proving
Items~(2)-(3) we can assume $F$ is locally-free: in fact if it is singular
replace it by $F^{**}$, and observe that stability for $F$ and $F^{**}$ is the
same, and that $|v(F^{**})|<|v(F)|$.  We will prove Item~(2) by contradiction.
Assume that the restriction  $F|_{C_t}$ to
the generic elliptic fiber $C_t$ is not stable. Since  $\la v^1(F),C\ra$ and
$v^0(F)$ are coprime the restriction $F|_{C_t}$ is unstable. Thus there
exists a proper subsheaf $A\ss F$ whose restriction to $C_t$ is a
desemistabilizing subsheaf, for the generic elliptic fiber $C_t$.
Hence $\left[r_F c_1(A)-r_A c_1(F)\ri]\cdot C>0$. If $N$ is sufficiently
large then $\left[r_F c_1(A)-r_A c_1(F)\ri]\cdot(NC+H)>0$, and hence $F$ is
$(NC+H)$-unstable. Since $H$ is $v$-stabilizing, $F$ is $H$-stable, and
hence by Lemma~(IV.3) the polarizations $H$ and
$(NC+H)$ are separated by a $|v|$-wall, contradicting the hypothesis that
$H$ is $|v|$-suitable. Thus $F|_{C_t}$ is stable for the generic elliptic fiber
$C_t$.  Now let's prove Item~(3). First let's show  that  $F$ is
$(NC+H)$-stable, for $N$ sufficiently large. In fact let $D\in |H|$ be a
smooth curve (if no such $D$ exists we replace $H$ by a high multiple of
itself), and set
$$\,\,d:=\max_{0\not=E\ss F}\{r_F\deg_D E-r_E\deg_D F\}\,,$$
where $\deg_D E=c_1(E|_D)$, $\deg_D F=c_1(F|_D)$.  Let
$C_1,\ldots,C_N \in |C|$ be generic elliptic fibers, and $A\ss F$ be a
subsheaf with $0<r_A<r_F$. Then
$$\,\,\left[r_Fc_1(A)-r_A c_1(F)\ri]\cdot(NC+H)=
r_F\deg_D A-r_A\deg_D F+
\sum_{i=1}^N r_F\deg_{C_i} A-r_A\deg_{C_i} F\,.$$
By stability of $F|_{C_i}$ the right-hand side is bounded above by $(d-N)$,
which is negative because  $N\gg 0$. Since $A\ss F$ is  arbitrary
we conclude that  $F$ is $(NC+H)$-stable. Proceeding as in the proof of
Item~(2) we conclude that $F$ is also $H$-stable.
This finishes the proof of Proposition~(I.1.6).
\bsk

\centerline{\bf References.}
\msk
\item {[A]} M.~F.~Atiyah. {\it Vector bundles over an elliptic curve},
Proc.~London Math.~Soc.~7 (1957), 414-452.
\item {[B]} A.~Beauville. {\it Vari\'et\'es K\"ahl\'eriennes dont la
premi\`ere classe de Chern est nulle}, J.~Differential Geom.~18 (1983),
755-782.
\item {[Fu]} A.~Fujiki. {\it On the de Rham cohomology group of compact
K\"ahler symplectic manifolds}, Algebraic Geometry Sendai 1985 (1987),
Adv. Studies in Pure Math. no. 10, Kinokuniya Tokio and North-Holland
Amsterdam, 105-165.
\item {[F1]} R.~Friedman. {\it Rank two vector bundles over regular
elliptic surfaces}, Invent.~math.~96 (1989), 283-332.
\item {[F2]} R.~Friedman. {\it Vector bundles and $SO(3)$-invariants for
elliptic surfaces III: the case of odd fiber degree}, preprint (1994).
\item {[Li1]} J.~Li. {\it The first two Betti numbers of the moduli spaces of
vector bundles on surfaces}, preprint (1995).
\item {[Li2]} J.~Li. {\it Picard groups of the moduli spaces of vector
bundles over algebraic surfaces}, preprint (1995).
\item{[GoHu]} L.~G\"ottsche, D.~Huybrechts. {\it Hodge numbers of moduli
spaces of stable bundles on $K3$ surfaces}, preprint, MPI/94-80.
\item{[GrHa]} P.~Griffiths, J.~Harris. {\it Principles of algebraic geometry},
John Wiles \& sons, 1978.

\item {[Ma]} M. Maruyama. {\it Moduli of stable sheaves II}, J. Math. Kyoto
Univ. 18 (1978), 557-614.
\item {[Mo]} J.~Morgan. {\it Comparison of the Donaldson polynomial
invariants with their algebro-geomteric analogues}, Topology 32 (1993),
449-489.
\item {[M1]} S.~Mukai. {\it Symplectic structure of the moduli
space of sheaves on an abelian or $K3$ surface}, Invent.~math.~77 (1984),
101-116.
\item {[M2]} S.~Mukai. {\it On the moduli space of bundles on $K3$ surfaces
I}, M. F. Atiyah and others, Vector bundles on algebraic varieties, Bombay
colloquium 1984, Tata Institute for fundamental research studies in
mathematics no.~11 (1987), 341-413.
\item {[M3]} S.~Mukai. {\it Moduli of vector bundles on $K3$ surfaces, and
symplectic manifolds}, Sugaku Expositions, vol.~1, n.~2 (1988), 139-174.
\item{[MR]} V.~B.~Mehta and A.~Ramanathan. {\it Restriction of stable
sheaves and representations of the fundamental group}, Invent.~Math.~77
(1984), 163-172.
\item {[O]} K.~O'Grady. {\it Relations among Donaldson polynomials of
certain algebraic surfaces}, to appear on Forum Mathematicum.
\bsk
\n
Universit\`a di Salerno, Facolt\`a di Scienze, Baronissi (Sa) - Italia
\msk
\n
{\tt ogrady@mat.uniroma1.it}

\end